\documentclass[a4paper,11pt]{article}
\usepackage{latexsym,amssymb,enumerate,amsmath,epsfig,amsthm,dsfont}
\usepackage[margin=1in]{geometry}
\usepackage{setspace,color}
\usepackage{tikz}
\usepackage{floatrow}
\usepackage{graphicx,subfigure}
\usepackage[ruled]{algorithm2e}
\usepackage{epstopdf}
\usepackage{grffile}
\newcommand{\x}{\mathbf{x}}
\newcommand{\y}{\mathbf{y}}
\newcommand{\z}{\mathbf{z}}
\newcommand{\p}{\mathbf{p}}
\newcommand{\A}{\mathbf{A}}
\newcommand{\bP}{\mathbf{P}}
\newcommand{\bS}{\mathbf{S}}
\newcommand{\q}{\mathbf{q}}

\newcommand{\D}{\mathbf{D}}
\newcommand{\I}{\mathbf{I}}
\newcommand{\RR}{\mathds{R}}
\newcommand{\mT}{\mathcal{T}}
\newcommand{\cof}{\mathrm{cof}}
\newtheorem{rmrk}{Remark}[section]
\begin{document}

\title{On the Numerical Solution of Nonlinear Eigenvalue Problems for the Monge-Amp\`{e}re Operator}%\thanks{Dedicated to Enrique Zuazua on the occasion of his 60th birthday.}% At most 5 thanks
\date{\today}
\author{Roland Glowinski \thanks{Department of Mathematics, University of Houston, 4800 Calhoun Road, Houston, TX77204, USA, and Department of Mathematics, the Hong Kong Baptist University, Hong Kong (Email: {\bf roland@math.uh.edu})}
\and
 Shingyu Leung\thanks{Department of Mathematics, the Hong Kong University of Science and Technology, Clear Water Bay, Hong Kong (Email: {\bf masyleung@ust.hk})}
 \and
 Hao Liu \thanks{School of Mathematics, Georgia Institute of Technology, 686 Cherry Street, Atlanta, GA 30332 USA (Email: {\bf hao.liu@math.gatech.edu})}
 \and
  Jianliang Qian \thanks{Department of Mathematics, Michigan State University, East Lansing, MI 48824 (Email: {\bf qian@math.msu.edu}).}\vspace{0.5cm}\\
  \emph{Dedicated to Enrique Zuazua on the occasion of his 60th birthday.}
}

\maketitle
\begin{abstract}
In this article, we report the results we obtained when investigating the numerical solution of some nonlinear eigenvalue problems for the Monge-Amp\`{e}re operator $v\rightarrow \det \D^2 v$. The methodology we employ relies on the following ingredients: (i) A divergence formulation of the eigenvalue problems under consideration. (ii) The time discretization by operator-splitting of an initial value problem (a kind of gradient flow) associated with each eigenvalue problem. (iii) A finite element approximation relying on spaces of continuous piecewise
affine functions. To validate the above methodology, we applied it to the solution of problems with known exact solutions: The results we obtained suggest convergence to the exact solution when the space discretization step $h\rightarrow 0$. We considered also test problems with no known exact solutions.

\end{abstract}

%%-----------------------------
%%      your text
%%-----------------------------
\section{Introduction}
\label{sec.intro}
There is an abundant literature concerning the solution of \emph{nonlinear eigenvalue problems} such as
\begin{align}
  -\nabla^2 u=\lambda f(u),
  \label{eq.ev}
\end{align}
where $\lambda$ is a real number and $f$ a real-valued function. These nonlinear problems include the celebrated \emph{Bratu-Gelfand} problem
\begin{align}
  -\nabla^2 u=\lambda e^u,
  \label{eq.bratu}
\end{align}
where $\lambda>0$. As shown in \cite{bebernes2013mathematical}, (\ref{eq.bratu}) has applications in \emph{solid combustion}. Problem (\ref{eq.bratu}) has motivated a very large number of publications, either mathematical or numerical (see, e.g., GOOGLE SCHOLAR for related references). On the other hand, to the best of our knowledge, there are only few publications on the Monge-Amp\`{e}re analogue of problem (\ref{eq.ev}), namely
\begin{align}
  \det \D^2u=\lambda f(u),
\end{align}
where $\D^2=\left( \frac{\partial^2}{\partial x_i \partial x_j}\right)_{1\leq i,j\leq d}$. Actually, the two publications we are aware of are \cite{lions1985two} and \cite{le2017eigenvalue}, where one discusses in particular the existence, uniqueness and regularity properties of the \emph{ground state} solutions to
\begin{align}
  \begin{cases}
    u\neq 0 \mbox{ and convex}, \lambda>0,\\
    \det\D^2 u=\lambda|u|^d \mbox{ in }\Omega,\\
    u=0 \mbox{ on }\partial\Omega,
  \end{cases}
  \label{eq.maevd}
\end{align}
where, in (\ref{eq.maevd}), $\Omega$ is a \emph{bounded convex domain} of $\RR^d$, with $d>1, \partial \Omega$ is the boundary of $\Omega$, $\lambda$ being the smallest nonlinear eigenvalue solution of problem (\ref{eq.maevd}). Problem (\ref{eq.maevd}) is degenerated since the right-hand side of the Monge-Amp\`{e}re type equation in (\ref{eq.maevd}) vanishes on $\partial\Omega$. Our goal in this publication is to address the \emph{numerical solution} of problem (\ref{eq.maevd}) and of the following variants of it:
\begin{align}
  \begin{cases}
    u\neq 0 \mbox{ and convex}, \lambda>0,\\
    \det\D^2 u=\lambda|u| \mbox{ in } \Omega,\\
    u=0 \mbox{ on } \partial \Omega,
  \end{cases}
  \label{eq.maev}
\end{align}
and
\begin{align}
  \begin{cases}
    u\neq 0 \mbox{ and convex}, \lambda>0,\\
    \det\D^2 u=\lambda e^{|u|} \mbox{ in } \Omega,\\
    u=0 \mbox{ on } \partial \Omega.
  \end{cases}
  \label{eq.mabg}
\end{align}
Clearly, problem (\ref{eq.mabg}) is related to the Bratu-Gelfand problem (\ref{eq.bratu}).

In Section \ref{sec.formulation}, we will provide \emph{divergence formulations} of problem (\ref{eq.maevd}), (\ref{eq.maev}) and (\ref{eq.mabg}), well-suited to their numerical solution by variational methods such as finite elements, and focusing on ground state solutions, reformulate these problems as minimization problems on nonlinear manifolds. In Section \ref{sec.timediscretize} we will discuss the time discretization by \emph{operator-splitting} of initial value problems associated with the above eigenvalue problems. In Section \ref{sec.spacediscretize}, we will discuss the \emph{mixed finite element approximation} of problems (\ref{eq.maevd}), (\ref{eq.maev}) and (\ref{eq.mabg}), with the Monge-Amp\`{e}re part of these problems treated by the methods discussed in \cite{glowinski2019finite} and \cite{liu2019finite}. Finally, in Section \ref{sec.numerical}, we will report on the results of \emph{numerical experiments}, validating the methodology discussed in the proceeding sections.

We dedicate this article to \emph{Professor Enrique Zuazua}, although (to the best of our knowledge) he has not contributed (yet) to the mathematics or the numerics of \emph{second order fully nonlinear elliptic equations}. Based on his outstanding scientific curiosity and capabilities, we have no doubt that Prof. E. Zuazua would significantly contribute to the above topics if for some reason they become (as we hope) of interest to him. Actually, the first author (RG) had the honor to have E. Zuazua attending the graduate course he was giving at Paris VI University in the eighties, and to witness, during the following years, his evolution from an excellent graduate student to a great scientist.
\section{Alterative formulations of problems (\ref{eq.maevd}), (\ref{eq.maev}) and (\ref{eq.mabg})}
\label{sec.formulation}

Denote by $\cof\left(\D^2 \phi\right)$ the cofactor matrix of $\D^2 \phi$. Using the identity
\begin{align}
  \nabla \cdot \left(\cof\left(\D^2 \phi\right)\nabla \phi\right)\equiv d\det \D^2\phi, \forall \phi \mbox{ smooth enough defined in }\RR^d,
\end{align}
and the fact that the convexity of $\Omega$ and $u$ and the condition $u|_{\partial\Omega}=0$ imply $u\leq 0$, we can reformulate problems (\ref{eq.maevd}), (\ref{eq.maev}) and (\ref{eq.mabg}) as
\begin{align}
  \begin{cases}
    u\leq 0 \mbox{ and convex}, \lambda>0,\\
    -\nabla \cdot \left(\cof\left(\D^2 u\right)\nabla u\right)=d\lambda u|u|^{d-1} \mbox{ in } \Omega,\\
    u=0 \mbox{ on } \partial \Omega,
  \end{cases}
  \label{eq.maevd.div}
\end{align}
\begin{align}
  \begin{cases}
    u\leq 0 \mbox{ and convex}, \lambda>0,\\
     -\nabla \cdot \left(\cof\left(\D^2 u\right)\nabla u\right)=d\lambda u \mbox{ in } \Omega,\\
    u=0 \mbox{ on } \partial \Omega,
  \end{cases}
  \label{eq.maev.div}
\end{align}
and
\begin{align}
  \begin{cases}
    u\leq 0 \mbox{ and convex}, \lambda>0,\\
     -\nabla \cdot \left(\cof\left(\D^2 u\right)\nabla u\right)=-d\lambda e^{-u} \mbox{ in } \Omega,\\
    u=0 \mbox{ on } \partial \Omega,
  \end{cases}
  \label{eq.mabg.div}
\end{align}
respectively. As in \cite{lions1985two} and \cite{le2017eigenvalue}, we are going to focus on \emph{ground state solutions} in the following way: in (\ref{eq.maevd.div}), (\ref{eq.maev.div}) and (\ref{eq.mabg.div}), we are going to consider $d\lambda$ as a \emph{Lagrange multiplier} associated with the following nonlinearly constrained problems from Calculus of Variations
\begin{subequations}
  \begin{align}
    u=\arg\min_{v\in S_1^d} \int_{\Omega} \left(\cof\left(\D^2 v\right)\nabla v\right)\cdot \nabla vd\x,
    \label{eq.maevd.min}
  \end{align}
  with $\x=\{x_i\}_{i=1}^d,d\x=dx_1\cdots dx_d$, and
  \begin{align}
    S_1^d=\left\{v\in V_d, v|_{\partial\Omega}=0, v\mbox{ convex}, \int_{\Omega} |v|^{d+1}d\x=C\right\},\ C\mbox{ being a positive constant,}
    \label{eq.maevd.s}
  \end{align}
  \label{eq.maevd.all}
\end{subequations}
then
\begin{subequations}
  \begin{align}
    u=\arg\min_{v\in S_2^d} \int_{\Omega} \left(\cof\left(\D^2 v\right)\nabla v\right)\cdot \nabla vd\x,
    \label{eq.maev.min}
  \end{align}
  with
  \begin{align}
    S_2^d=\left\{v\in V_d, v|_{\partial\Omega}=0, v\mbox{ convex}, \int_{\Omega} |v|^2d\x=C\right\},\ C\mbox{ being a positive constant,}
    \label{eq.maev.s}
  \end{align}
  \label{eq.maev.all}
\end{subequations}
and
\begin{subequations}
  \begin{align}
    u=\arg\min_{v\in S_3^d} \int_{\Omega} \left(\cof\left(\D^2 v\right)\nabla v\right)\cdot \nabla vd\x,
    \label{eq.mabg.min}
  \end{align}
  with
  \begin{align}
    S_3^d=\left\{v\in V_d, v|_{\partial\Omega}=0, v\mbox{ convex}, \int_{\Omega} \left(e^{-v}-1\right)d\x=C\right\},\ C\mbox{ being a positive constant,}
    \label{eq.mabg.s}
  \end{align}
  \label{eq.mabg.all}
\end{subequations}
respectively. From \emph{Sobolev imbedding theorems} (see, e.g., Chapter 6 of \cite{ciarlet2013linear}, and the references therein) the 'largest' and simplest \emph{Sobolev spaces} $V_d$, for which the various integrals in (\ref{eq.maevd.all}), (\ref{eq.maev.all}) and (\ref{eq.mabg.all}) make sense are
\begin{align}
  \begin{cases}
    V_2=W^{2,\frac{3}{2}}(\Omega),\\
    V_3=W^{2,\frac{12}{5}}(\Omega).
  \end{cases}
  \label{eq.V}
\end{align}
If the spaces $V_d$ are defined by (\ref{eq.V}), the integrands in (\ref{eq.maevd.min}), (\ref{eq.maev.min}) and (\ref{eq.mabg.min}) belong to $L^1(\Omega)$. Several remarks are in order; among them
\begin{rmrk}
  From a mathematical point of view, the transition from (\ref{eq.maevd.div}), (\ref{eq.maev.div}) and (\ref{eq.mabg.div}) to (\ref{eq.maevd.all}), (\ref{eq.maev.all}) and (\ref{eq.mabg.all}) is formal and deserves a mathematical justification (not given here). On the other hand, the numerical results reported in Section \ref{sec.numerical} validate employing formulations (\ref{eq.maevd.all}), (\ref{eq.maev.all}) and (\ref{eq.mabg.all}) to compute the ground state solutions of problems (\ref{eq.maevd.div}), (\ref{eq.maev.div}) and (\ref{eq.mabg.div}).
\end{rmrk}

\begin{rmrk}
  Suppose that the pair $(u,\lambda)$ is solution to problem (\ref{eq.maevd}), (\ref{eq.maevd.div}), then the pair $(\theta u,\lambda)$ is also solution to (\ref{eq.maevd}), (\ref{eq.maevd.div}), for any positive number $\theta$. Similarly, if the pair $(u,\lambda)$ is solution to problem (\ref{eq.maev}), (\ref{eq.maev.div}), the function $w=\theta u$ verifies
  \begin{align*}
  \begin{cases}
    w\neq 0 \mbox{ and convex}, \lambda>0,\\
    \det \D^2 w=\theta^{d-1}\lambda |w| \mbox{ in }\Omega, \forall \theta>0,\\
    w=0 \mbox{ on } \partial\Omega,
  \end{cases}
  \end{align*}
  implying that the pair $(w,\theta^{d-1}\lambda)$ is solution of a nonlinear eigenvalue problem of type (\ref{eq.maev}), (\ref{eq.maev.div}). These observations suggest taking $C=1$ in (\ref{eq.maevd.s}) and (\ref{eq.maev.s}).
\end{rmrk}

\begin{rmrk}
  Suppose that $\Omega=\left\{ (x_1,x_2)\in \RR^2,x_1^2+x_2^2<1\right\}$. The \emph{radial} solutions of problem (\ref{eq.mabg}) verify
  \begin{align}
    \begin{cases}
      u\leq0 \mbox{ on } [0,1], \lambda>0,\\
      u'u''=\lambda r e^{-u} \mbox{ on } (0,1),\\
      u(0)<0, u'(0)=0, u(1)=0,
    \end{cases}
    \label{eq.mabg.radial}
  \end{align}
  with $r=\sqrt{x_1^2+x_2^2}$. A simple way to solve (\ref{eq.mabg.radial}) is to pick $u(0)$ and to use a shooting method where one adjusts $\lambda$ iteratively in order to obtain $u(1)=0$. Doing so, we obtain the \emph{bifurcation diagram} visualized in Figure \ref{fig.mabg.bifurcation}, a diagram qualitatively similar to those associated with the Bratu problem
  $$
  -\nabla^2 u =\lambda e^u \mbox{ in }\Omega, u=0 \mbox{ on } \partial\Omega, \lambda>0,
  $$
  if $\Omega$ is a bounded convex domain of $\RR^2$. We will take advantage of these results to validate the methodology discussed in the following sections.
  \begin{figure}[t]
  \includegraphics[width=0.6\textwidth]{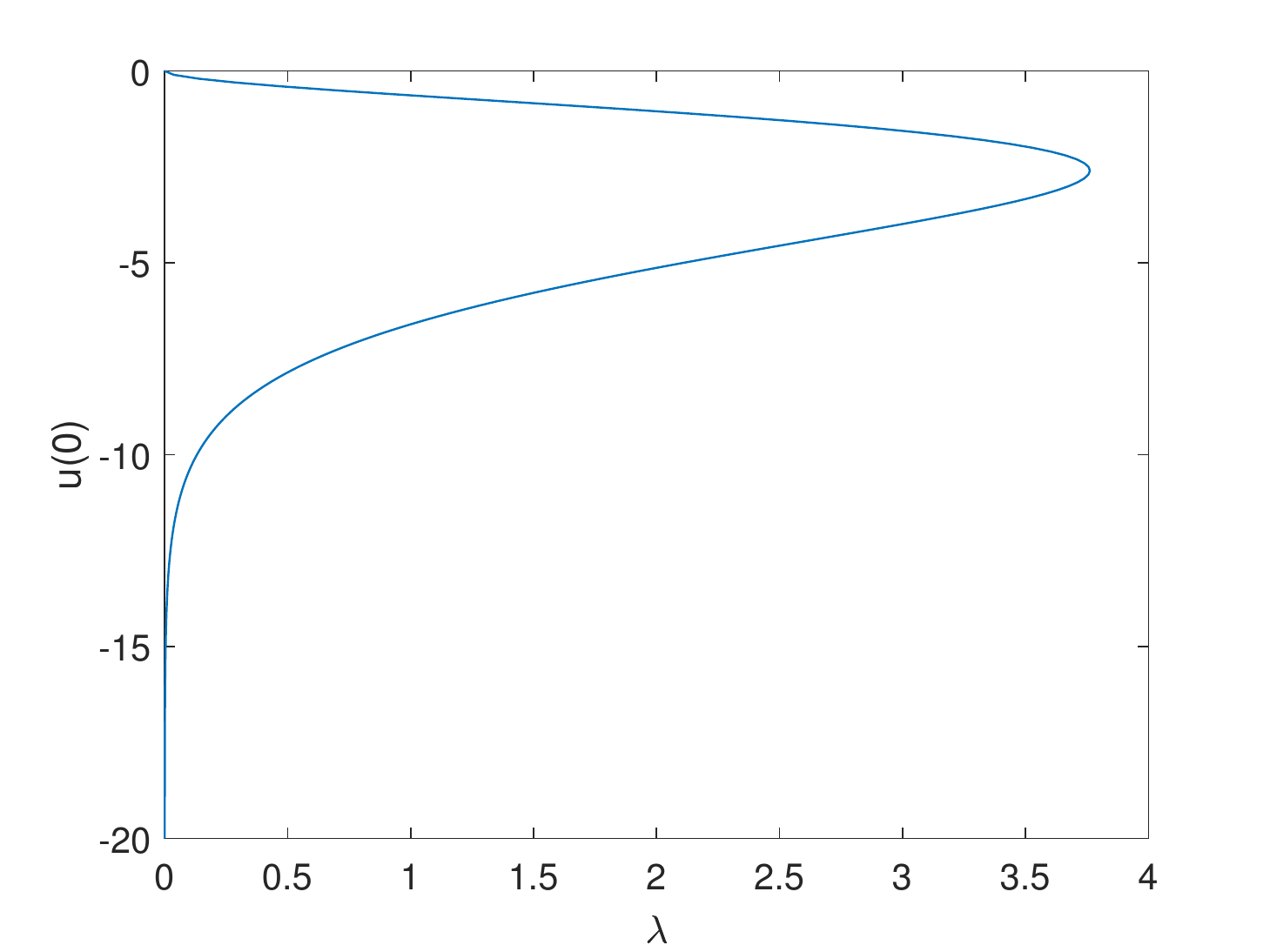}
    \caption{Problem (\ref{eq.mabg.radial}) bifurcation diagram: We have visualized the variations of $u(0)$ (the minimal value of $u$) as a function of $\lambda$. The turning point is located at $\lambda\approx 3.7617$ and $u(0)\approx-2.5950$.}
    \label{fig.mabg.bifurcation}
  \end{figure}
\end{rmrk}

Written in PDE form, the problems we are going to address amount to finding the ground state solutions of
\begin{align}
  \begin{cases}
    u\leq 0 \mbox{ and convex}, \lambda>0,\\
    \begin{cases}
      -\nabla \cdot \left(\cof\left(\D^2 u\right)\nabla u\right)=d\lambda u|u|^{d-1} \mbox{ in } \Omega,\\
      u=0 \mbox{ on } \partial\Omega,
    \end{cases}\\
    {\displaystyle \int_{\Omega}}|u|^{d+1}d\x=1,
  \end{cases}
  \label{eq.maevd.sys}
\end{align}
\begin{align}
  \begin{cases}
    u\leq 0 \mbox{ and convex}, \lambda>0,\\
    \begin{cases}
      -\nabla \cdot \left(\cof\left(\D^2 u\right)\nabla u\right)=d\lambda u \mbox{ in } \Omega,\\
      u=0 \mbox{ on } \partial\Omega,
    \end{cases}\\
    {\displaystyle \int_{\Omega}}|u|^2d\x=1,
  \end{cases}
  \label{eq.maev.sys}
\end{align}
and
\begin{align}
  \begin{cases}
    u\leq 0 \mbox{ and convex}, \lambda>0,\\
    \begin{cases}
      -\nabla \cdot \left(\cof\left(\D^2 u\right)\nabla u\right)=-d\lambda e^{-u} \mbox{ in } \Omega,\\
      u=0 \mbox{ on } \partial\Omega,
    \end{cases}\\
    {\displaystyle \int_{\Omega}}\left(e^{-u}-1\right)d\x=C(>0).
  \end{cases}
  \label{eq.mabg.sys}
\end{align}
\section{On an operator-splitting method for the solution of problems (\ref{eq.maevd.sys}), (\ref{eq.maev.sys}) and (\ref{eq.mabg.sys})}
\label{sec.timediscretize}
\subsection{On the solution of a simple linear eigenvalue problem}
\label{sec.timediscretize.linearEV}
In order to motivate the approach we are going to employ, let us consider the following \emph{linear eigenvalue problem}:
\begin{align}
  \lambda=\min_{\y\in S} \A\y\cdot \y,
  \label{eq.linearEV}
\end{align}
where in (\ref{eq.linearEV}): (i) $\A$ is a $d\times d$ real symmetric matrix. (ii) $\y\cdot \z=\sum_{i=1}^d y_iz_i,\forall \y=(y_i)_{i=1}^d,\z=(z_i)_{i=1}^d\in \RR^d$. (iii) $S=\left\{ \y\in \RR^d, \|\y\|_2=\sqrt{\y\cdot \y}=1\right\}$, implying that $\lambda$ is the \emph{smallest eigenvalue} of matrix $\A$. If $\x\in S$ is an eigenvector of $\A$ associated with $\lambda$, $\x$ is a soluton of the following nonlinearly constrained minimization problem
\begin{align}
  \x\in \arg\min_{\y\in S} \frac{1}{2}\A \y\cdot\y.
  \label{eq.linearEV.min}
\end{align}
Let $I_S$ be the \emph{indicator functional} of sphere $S$, that is
\begin{align}
  I_S(\y)=\begin{cases}
    0, &\mbox{ if } \y\in S,\\
    +\infty, &\mbox{ if } \y\in \RR^d\backslash S.
  \end{cases}
  \label{eq.indicator}
\end{align}
It follows then from (\ref{eq.linearEV.min}) and (\ref{eq.indicator}) that
\begin{align}
  \x\in \arg\min_{\y\in \RR^d} \left[ \frac{1}{2} \A\y\cdot \y+ I_S(\y)\right].
\end{align}
Formally, $\x$ verifies the following \emph{necessary optimality condition}
\begin{align}
  \A\x+\partial I_S(\x)\ni \mathbf{0},
  \label{eq.optimality}
\end{align}
where $\partial I_S(\x)$ denotes a \emph{generalized differential} of functional $I_S$ at $\x$ (for the notions of generalized differentials see, e.g., \cite{dutta2005generalized} and the references therein). In order to solve (\ref{eq.optimality}), we associate with it the following \emph{initial value problem} (\emph{gradient flow} in the dynamical system terminology):
\begin{align}
  \begin{cases}
    \x(t)\in \RR^d,t>0,\\
    \frac{d\x}{dt}+\A\x+\partial I_S(\x)\ni \mathbf{0} \mbox{ on } (0,+\infty),\\
    \x(0)=\x_0(\in\RR^d).
  \end{cases}
  \label{eq.linearEV.ivp}
\end{align}
Following, e.g., \cite{glowinski2017splitting}, we apply to (\ref{eq.linearEV.ivp}) the \emph{Marchuk-Yanenko scheme} (with $\tau(>0)$ a time-discretization step and $\x^n$ an approximation of $\x(n\tau)$), that is
\begin{align}
  \x^0=\x_0.
  \label{eq.linearEV.split.0}
\end{align}
For $n\geq 0,\x^n\rightarrow \x^{n+1/2}\rightarrow \x^{n+1}$ as follows:
\begin{align}
  \frac{\x^{n+1/2}-\x^n}{\tau}+\A\x^{n+1/2}=\mathbf{0},
  \label{eq.linearEV.split.1}
\end{align}
and
\begin{align}
  \frac{\x^{n+1}-\x^{n+1/2}}{\tau}+\partial I_S(\x^{n+1})\ni \mathbf{0}.
  \label{eq.linearEV.split.2}
\end{align}
We clearly have
\begin{align}
  \x^{n+1/2}=(\I+\tau\A)^{-1}\x^n
  \label{eq.linearEV.split.solve1}
\end{align}
(assuming, of course, that $\det(\I+\tau\A)\neq 0$).

From the context of the problem, we interpret (\ref{eq.linearEV.split.2}) as a necessary optimality condition for the following minimization problem:
\begin{align}
  \x^{n+1}\in \arg\min_{\y\in S} \left[ \frac{1}{2} \y\cdot\y-\x^{n+1/2}\cdot \y\right].
  \label{eq.linearEV.split.solve2}
\end{align}
Since $\y\cdot\y=1$ over $S$, problem (\ref{eq.linearEV.split.solve2}) reduces to
\begin{align}
  \x^{n+1}\in \arg\max_{\y\in S} \x^{n+1/2}\cdot \y,
\end{align}
whose solution is
\begin{align}
  \x^{n+1}=\begin{cases}
    \frac{\x^{n+1/2}}{\|\x^{n+1/2}\|_2} &\mbox{ if } \x^{n+1/2}\neq \mathbf{0},\\
    \mbox{any element of } S &\mbox{ if }\x^{n+1/2}=\mathbf{0}.
  \end{cases}
  \label{eq.linearEV.split.solve3}
\end{align}
Actually, $\x_0\neq 0\Rightarrow \x^{n+1/2}\neq 0, \forall n>0$. It follows from relations (\ref{eq.linearEV.split.solve1}) and (\ref{eq.linearEV.split.solve3}) that algorithm (\ref{eq.linearEV.split.0})-(\ref{eq.linearEV.split.2}) is a disguised form of the \emph{inverse power method with shift} for \emph{eigenvalue} and \emph{eigenvector computation}, a well-known method from \emph{numerical linear algebra}. Let us denote by $\lambda_1(=\lambda_1^+-\lambda_1^-)$ the \emph{smallest eigenvalue} of matrix $\A$ and by $\pi_1(\x_0)$ the orthogonal projection of $\x_0$ over the subspace of $\RR^d$ spanned by the eigenvectors of $\A$ associated with $\lambda_1$. Assume that $\pi_1(\x_0)\neq \mathbf{0}$; we have then the following convergence result:
\begin{align}
  \tau<1/\lambda_1^-\Rightarrow \lim_{n\rightarrow +\infty} \x^n=\x=\frac{\pi_1(\x_0)}{\|\pi_1(\x_0)\|_2}.
\end{align}
The algorithms we are going to use, to compute the ground state solutions of the nonlinear eigenvalue problems described in Section \ref{sec.formulation}, are conceptually close to algorithm (\ref{eq.linearEV.split.0})-(\ref{eq.linearEV.split.2}).

\subsection{Initial value problems associated with the nonlinear eigenvalue problems (\ref{eq.maevd.sys}), (\ref{eq.maev.sys}) and (\ref{eq.mabg.sys})}
To facilitate the solution of problems (\ref{eq.maevd.sys}), (\ref{eq.maev.sys}) and (\ref{eq.mabg.sys}) via a time dependent approach inspired by Section \ref{sec.timediscretize.linearEV}, we reformulate the above three problems as:
\begin{align}
  \begin{cases}
    u\leq 0, \p \mbox{ pointwise positive semi-definite over } \Omega, \lambda>0,\\
    \begin{cases}
      -\nabla \cdot \left(\cof\left(\p\right)\nabla u\right)=d\lambda u|u|^{d-1} \mbox{ in } \Omega,\\
      u=0 \mbox{ on } \partial\Omega,
    \end{cases}\\
    \p=\D^2 u,\\
    {\displaystyle \int_{\Omega}}|u|^{d+1}d\x=1,
  \end{cases}
  \label{eq.maevd.sys.1}
\end{align}
\begin{align}
  \begin{cases}
    u\leq 0, \p \mbox{ pointwise positive semi-definite over } \Omega, \lambda>0,\\
    \begin{cases}
      -\nabla \cdot \left(\cof\left(\p\right)\nabla u\right)=d\lambda u \mbox{ in } \Omega,\\
      u=0 \mbox{ on } \partial\Omega,
    \end{cases}\\
    \p=\D^2 u,\\
    {\displaystyle \int_{\Omega}}|u|^2d\x=1,
  \end{cases}
  \label{eq.maev.sys.1}
\end{align}
and
\begin{align}
  \begin{cases}
    u\leq 0, \p \mbox{ pointwise positive semi-definite over } \Omega, \lambda>0,\\
    \begin{cases}
      -\nabla \cdot \left(\cof\left(\p\right)\nabla u\right)=-d\lambda e^{-u} \mbox{ in } \Omega,\\
      u=0 \mbox{ on } \partial\Omega,
    \end{cases}\\
    \p=\D^2 u,\\
    {\displaystyle \int_{\Omega}}\left(e^{-u}-1\right)d\x=C(>0),
  \end{cases}
  \label{eq.mabg.sys.1}
\end{align}
respectively. We associate with (\ref{eq.maevd.sys.1}), (\ref{eq.maev.sys.1}) and (\ref{eq.mabg.sys.1}) the following initial value problems:
\begin{align}
&\mbox{Find } u(t)\leq0, \p(t) \mbox{ pointwise symmetric positive semi-definite, and } \lambda(t)>0, \mbox{ so that} \nonumber\\
 &\begin{cases}
    \begin{cases}
      \frac{\partial u}{\partial t}-\nabla\cdot \left[\left(\varepsilon \I+\cof(\p)\right)\nabla u\right] =d\lambda u|u|^{d-1} \mbox{ in } \Omega\times (0,+\infty),\\
      u=0 \mbox{ on } \partial\Omega\times (0,+\infty),
    \end{cases}\\
    \frac{\partial \p}{\partial t}+\gamma(\p-\D^2 u)=\mathbf{0} \mbox{ in } \Omega\times (0,+\infty),\\
    {\displaystyle \int_{\Omega}} |u(t)|^{d+1} d\x =1, \forall t>0,\\
    (u(0),\p(0))=(u_0,\p_0),
  \end{cases}
  \label{eq.maevd.ivp}
\end{align}
\begin{align}
&\mbox{Find } u(t)\leq0, \p(t) \mbox{ pointwise symmetric positive semi-definite, and } \lambda(t)>0, \mbox{ so that} \nonumber\\
 &\begin{cases}
    \begin{cases}
      \frac{\partial u}{\partial t}-\nabla\cdot \left[\left(\varepsilon \I+\cof(\p)\right)\nabla u\right] =d\lambda u \mbox{ in } \Omega\times (0,+\infty),\\
      u=0 \mbox{ on } \partial\Omega\times (0,+\infty),
    \end{cases}\\
    \frac{\partial \p}{\partial t}+\gamma(\p-\D^2 u)=\mathbf{0} \mbox{ in } \Omega\times (0,+\infty),\\
    {\displaystyle \int_{\Omega}} |u(t)|^2 d\x =1, \forall t>0,\\
    (u(0),\p(0))=(u_0,\p_0),
  \end{cases}
  \label{eq.maev.ivp}
\end{align}
and
\begin{align}
&\mbox{Find } u(t)\leq0, \p(t) \mbox{ pointwise symmetric positive semi-definite, and } \lambda(t)>0, \mbox{ so that} \nonumber\\
 &\begin{cases}
    \begin{cases}
      \frac{\partial u}{\partial t}-\nabla\cdot \left[\left(\varepsilon \I+\cof(\p)\right)\nabla u\right] =-d\lambda e^{-u} \mbox{ in } \Omega\times (0,+\infty),\\
      u=0 \mbox{ on } \partial\Omega\times (0,+\infty),
    \end{cases}\\
    \frac{\partial \p}{\partial t}+\gamma(\p-\D^2 u)=\mathbf{0} \mbox{ in } \Omega\times (0,+\infty),\\
    {\displaystyle \int_{\Omega}} \left(e^{-u(t)}-1\right) d\x =C(>0), \forall t>0,\\
    (u(0),\p(0))=(u_0,\p_0),
  \end{cases}
  \label{eq.mabg.ivp}
\end{align}
respectively. In (\ref{eq.maevd.ivp})-(\ref{eq.mabg.ivp}): (i) $\phi(t)$ denotes the function $\x\rightarrow \phi(\x,t)$. (ii) $\varepsilon>0$ and is of the order of $h^2$ in practice ($h$ being a space discretization step). (iii) $\gamma\geq \beta\lambda_0\varepsilon$, with $\beta$ a positive number of the order of $1$ and $\lambda_0(>0)$ the \emph{smallest eigenvalue} of operator $-\nabla^2$ over the space $H_0^1(\Omega)$ (see \cite{glowinski2019finite} for the rational of this choice of $\gamma$). (iv) $u_0\leq 0$ and $\neq 0$. (v) $\p^0$ is pointwise symmetric positive semi-definite.

Problem (\ref{eq.maev.ivp}) being the simplest of the three above initial value problems, it will be the first whose operator-splitting solution will be discussed (in Section \ref{sec.timediscretize.maev}).

\subsection{Operator-splitting solution of problem (\ref{eq.maev.ivp})}
\label{sec.timediscretize.maev}
With $\tau(>0)$ a \emph{time-discretization step} (fixed for simplicity) we approximate the initial value problem (\ref{eq.maev.ivp}) by
\begin{align}
  (u^0,\p^0)=(u_0,\p_0).
  \label{eq.maev.split.0}
\end{align}
For $n\geq0, (u^n,\p^n)\rightarrow u^{n+1/2}\rightarrow (u^{n+1},\p^{n+1})$ as follows:\\
\underline{\emph{First Step}}: Solve the following \emph{(well-posed) linear elliptic problem}
\begin{align}
  \begin{cases}
    u^{n+1/2}-\tau \nabla\cdot \left[ \left(\varepsilon\I+\cof(\p^n)\right)\nabla u^{n+1/2}\right] =u^n \mbox{ in } \Omega,\\
    u^{n+1/2}=0 \mbox{ on }\partial\Omega.
  \end{cases}
  \label{eq.maev.split.1}
\end{align}
\underline{\emph{Second Step}}:
\begin{subequations}
  \begin{align}
    &\p^{n+1}(\x)=P_+\left[ e^{-\gamma\tau} \p^n(\x)+\left(1-e^{-\gamma\tau}\right)\D^2 u^{n+1/2}(\x)\right], \mbox{ a.e. } \x\in\Omega, \label{eq.maev.split.2a}\\
    &\begin{cases}
      u^{n+1}-u^{n+1/2}=\tau d\lambda^{n+1} u^{n+1},\\
      \displaystyle\int_{\Omega} |u^{n+1}|^2d\x=1\ \left(\Leftrightarrow u^{n+1}\in \Sigma_2=\left\{ \phi\in L^2(\Omega), \displaystyle\int_{\Omega}|\phi|^2d\x=1\right\}\right).
    \end{cases}
    \label{eq.maev.split.2b}
  \end{align}
  \label{eq.maev.split.2}
\end{subequations}
Above:
\begin{itemize}
  \item Problem (\ref{eq.maev.split.1}) is a classical \emph{linear elliptic problem}. Its unique solution $u^{n+1/2}$ has the (easy to prove) following properties
      \begin{align}
        \begin{cases}
          u^n\neq 0,\\
          u^n\leq 0,
        \end{cases}\Rightarrow
        \begin{cases}
          u^{n+1/2}\neq 0,\\
          u^{n+1/2}\leq 0,
        \end{cases}
        \left\| u^{n+1/2}\right\|_{L^2(\Omega)} < \left\| u^n\right\|_{L^2(\Omega)}(=1 \mbox{ if }n\geq 1).
        \label{eq.maev.split.1.prop}
      \end{align}
  \item $P_+$ is an operator, mapping the space of the $d\times d$ \emph{real symmetric matrices} onto the \emph{closed convex cone of the real symmetric positive semi-definite $d\times d$ matrices}. If $\q$ is a $2\times 2$ real symmetric matrix with eigenvalues $\mu_1$ and $\mu_2$, $\exists \bS\in O(2)$ such that $\q=\bS\begin{pmatrix}\mu_1 & 0\\ 0 & \mu_2\end{pmatrix} \bS^{-1}$, then, operator $P_+$ is defined by
      \begin{align}
        P_+(\q)=\bS\begin{pmatrix}\max(0,\mu_1) & 0\\ 0 & \max(0,\mu_2)\end{pmatrix} \bS^{-1}.
      \end{align}
      One should proceed similarly if $d\geq 3$.
  \item Since we are looking for \emph{ground state solutions}, we will consider (\ref{eq.maev.split.2b}) as an \emph{optimality system} associated with the following \emph{minimization problem}
      \begin{align}
        u^{n+1}=\arg\min_{v\in \Sigma_2} \left[ \frac{1}{2} \int_{\Omega} |v|^2d\x-\int_{\Omega} u^{n+1/2} vd\x\right].
        \label{eq.maev.split.2b.solve}
      \end{align}
      Since $\displaystyle\int_{\Omega} |v|^2d\x=1$ over $\Sigma_2$, it follows from (\ref{eq.maev.split.2b.solve}) that
      \begin{align*}
        u^{n+1}=\arg\max_{v\in\Sigma_2} \int_{\Omega} u^{n+1/2}vd\x,
      \end{align*}
      that is (since, from (\ref{eq.maev.split.1.prop}), $u^{n+1/2}\neq 0$):
      \begin{align}
        u^{n+1}=\frac{u^{n+1/2}}{\left\| u^{n+1/2}\right\|_{L^2(\Omega)}}.
        \label{eq.maev.split.2b.solve2}
      \end{align}
      It follows from (\ref{eq.maev.split.2b}), (\ref{eq.maev.split.1.prop}) and (\ref{eq.maev.split.2b.solve2}) that
      \begin{align}
        \lambda^{n+1}=\frac{1-\left\|u^{n+1/2}\right\|_{L^2(\Omega)}}{d\tau}>0, \forall n\geq 1.
        \label{eq.maev.split.lambda}
      \end{align}
\end{itemize}
Several comments and remarks are in order concerning algorithm (\ref{eq.maev.split.0})-(\ref{eq.maev.split.2}). Among them:
\begin{rmrk}
  Algorithm (\ref{eq.maev.split.0})-(\ref{eq.maev.split.2}) has clearly the flavor of an \emph{inverse power method}.
\end{rmrk}
\begin{rmrk}
\label{remark.spliterror}
  Algorithm (\ref{eq.maev.split.0})-(\ref{eq.maev.split.2}) 'enjoys' a splitting error forcing us to use a small time-discretization step $\tau$ (actually, a stability condition associated with part (\ref{eq.maev.split.1}) of the scheme requires to take $\tau$ of the order of $h^2$; see \cite{glowinski2019finite} for details).
\end{rmrk}
\begin{rmrk}
\label{remark.lambda}
  There is no need to compute $\lambda^{n+1}$ at each time step to obtain $\lambda$. Indeed, it follows from (\ref{eq.maev.split.lambda}) that $\lambda^{n+1}$ is obtained by the \emph{ratio of two small numbers}. It is safer (we think so at least) to proceed as follows: Denote by $(u_{\tau},\p_{\tau})$ the limit of $(u^n,\p^n)$. It makes sense to approximate the (nonlinear) eigenvalue $\lambda$ by
  \begin{align*}
    \lambda_{\tau}=\frac{\displaystyle\int_{\Omega} (\varepsilon\I+\cof(\p_{\tau}))\nabla u_{\tau} \cdot \nabla u_{\tau}d\x}{d},
  \end{align*}
  a (kind of) generalized \emph{Rayleigh quotient} (other approximations of the same type are available).
\end{rmrk}
\begin{rmrk}
\label{remark.initial}
  Concerning the choice of $(u_0,\p_0)$ we suggest the following strategy
  \begin{enumerate}[(i)]
    \item Compute the \emph{convex solution} of
    $$
    \begin{cases}
      \det\D^2 \psi_0=1 \mbox{ in }\Omega,\\

      \psi_0=0 \mbox{ on } \partial\Omega,
    \end{cases}
    $$
    using, for example, the methods discussed in \cite{glowinski2019finite} and \cite{liu2019finite}. We have $\psi_0\leq 0$ and $\psi_0\neq 0$.
    \item Define the pair $(u_0,\p_0)$ by
    \begin{align}
      u_0=\frac{\psi_0}{\|\psi_0\|_{L^2(\Omega)}} \mbox{ and } \p_0=D^2 u_0.
      \label{eq.maev.initial}
    \end{align}
    By construction, one has $u_0\leq 0, u_0\in \Sigma_2$ and $\p_0$ \emph{symmetric positive definite}.
  \end{enumerate}
\end{rmrk}

\subsection{Operator-splitting solution of problem (\ref{eq.maevd.ivp})}
\label{sec.timediscretize.maevd}
We use the notation of Section \ref{sec.timediscretize.maev}. The variant of algorithm (\ref{eq.maev.split.0})-(\ref{eq.maev.split.2}) reads as:
\begin{align}
  (u^0,\p^0)=(u_0,\p_0).
  \label{eq.maevd.split.0}
\end{align}
For $n\geq0, (u^n,\p^n)\rightarrow u^{n+1/3}\rightarrow (u^{n+2/3},\p^{n+1})\rightarrow u^{n+1}$ as follows:\\
\underline{\emph{First Step}}: Solve the following \emph{(well-posed) linear elliptic problem}
\begin{align}
  \begin{cases}
    u^{n+1/3}-\tau \nabla\cdot \left[ \left(\varepsilon\I+\cof(\p^n)\right)\nabla u^{n+1/3}\right] =u^n \mbox{ in } \Omega,\\
    u^{n+1/3}=0 \mbox{ on }\partial\Omega.
  \end{cases}
  \label{eq.maevd.split.1}
\end{align}
\underline{\emph{Second Step}}:
\begin{subequations}
  \begin{align}
    &\p^{n+1}(\x)=P_+\left[ e^{-\gamma\tau} \p^n(\x)+\left(1-e^{-\gamma\tau}\right)\D^2 u^{n+1/3}(\x)\right], \mbox{ a.e. } \x\in\Omega, \label{eq.maev.split.2a}\\
    &\begin{cases}
      u^{n+2/3}-u^{n+1/3}=\tau d\lambda^{n+1} u^{n+2/3}|u^{n+2/3}|^{d-1},\\
      {\displaystyle \int_{\Omega}} |u^{n+2/3}|^{d+1}d\x=1\ \left(\Leftrightarrow u^{n+2/3}\in \Sigma_{d+1}=\left\{ \phi\in L^{d+1}(\Omega), {\displaystyle \int_{\Omega}}|\phi|^{d+1}d\x=1\right\}\right).
    \end{cases}
    \label{eq.maevd.split.2b}
  \end{align}
  \label{eq.maevd.split.2}
\end{subequations}
\underline{\emph{Third Step}}:
\begin{align}
  u^{n+1}=-|u^{n+2/3}|.
  \label{eq.maevd.split.3}
\end{align}
Above:
\begin{itemize}
  \item The following variant of relations (\ref{eq.maev.split.1.prop}) holds:
  \begin{align}
        \begin{cases}
          u^n\neq 0,\\
          u^n\leq 0,
        \end{cases}\Rightarrow
        \begin{cases}
          u^{n+1/3}\neq 0,\\
          u^{n+1/3}\leq 0,
        \end{cases}
        \left\| u^{n+1/3}\right\|_{L^{2}(\Omega)} < \left\| u^n\right\|_{L^{2}(\Omega)}.
        \label{eq.maevd.split.1.prop}
      \end{align}
  \item Since we are looking for \emph{ground state solutions}, we will consider (\ref{eq.maevd.split.2b}) as an \emph{optimality system} associated with the following \emph{minimization problem}
      \begin{align}
        u^{n+2/3}\in \arg\min_{v\in \Sigma_{d+1}} \left[ \frac{1}{2} \int_{\Omega} |v|^2d\x -\int_{\Omega} u^{n+1/3}vd\x\right].
        \label{eq.maevd.split.2.solve}
      \end{align}
      We are facing a mathematical difficulty with problem (\ref{eq.maevd.split.2.solve}): Indeed problem (\ref{eq.maevd.split.2.solve}) is \emph{ill-posed} since the sphere $\Sigma_{d+1}$ is not compact in $L^2(\Omega)$ if $d\geq 2$ (neither strongly, nor weakly). Fortunately, the \emph{finite element analogues} of (\ref{eq.maevd.split.2.solve}) we will encounter in Section \ref{sec.spacediscretize} have solutions. If problem (\ref{eq.maevd.split.2.solve}) has a solution $w$, one can easily show that $-|w|$ is also a solution, justifying relation (\ref{eq.maevd.split.3}). We are going to discuss now the (formal) solution of problem (\ref{eq.maevd.split.2.solve}) by a \emph{sequential quadratic programming} (SQP) algorithm, obtained by successive linearization of the constraint $\displaystyle\int_{\Omega} |v|^{d+1}d\x=1$ in the neighborhood of $u_k^{n+2/3}$ for $k\geq 0$. This SQP algorithm reads as:
      \begin{align}
        u_{0}^{n+2/3} =\begin{cases}
          u^0 &\mbox{ if } n=0,\\
          u^{(n-1)+2/3} &\mbox{ if } n\geq 1.
          \label{eq.maevd.SQP.0}
        \end{cases}
      \end{align}
      For $k\geq 0, u_k^{n+2/3}\rightarrow u_{k+1}^{n+2/3}$ as follows:
      \begin{align}
        &u_{k+1}^{n+2/3}=u^{n+1/3}+\nonumber\\
        &\left[\frac{1-\left\|u_k^{n+2/3}\right\|_{L^{d+1}(\Omega)}^{d+1}+ (d+1){\displaystyle\int_{\Omega}} \left( u_k^{n+2/3}-u^{n+1/3}\right) u_k^{n+2/3} \left|u_k^{n+2/3}\right|^{d-1} d\x}{(d+1)\left\|u_k^{n+2/3}\right\|_{L^{2d}(\Omega)}^{2d}}\right] u_k^{n+2/3} \left| u_k^{n+2/3}\right|^{d-1}.
        \label{eq.maevd.SQP.1}
      \end{align}
      Several variants of algorithm (\ref{eq.maevd.SQP.0})-(\ref{eq.maevd.SQP.1}) make sense: For example, since $\tau$ has to be small, one can stop iterating after just one iteration, implying that $u^{n+2/3}$ is defined from $u^{n+1/3}$ by
      \begin{align}
        &u^{n+2/3}=u^{n+1/3}+\nonumber\\
        &\left[\frac{1-\left\|u_0^{n+2/3}\right\|_{L^{d+1}(\Omega)}^{d+1}+ (d+1){\displaystyle\int_{\Omega}} \left( u_0^{n+2/3}-u^{n+1/3}\right) u_0^{n+2/3} \left|u_0^{n+2/3}\right|^{d-1} d\x}{(d+1)\left\|u_0^{n+2/3}\right\|_{L^{2d}(\Omega)}^{2d}}\right] u_0^{n+2/3} \left| u_0^{n+2/3}\right|^{d-1}.
      \end{align}
\end{itemize}

One obtains another variant of algorithm (\ref{eq.maevd.SQP.0})-(\ref{eq.maevd.SQP.1}) by introducing a \emph{damping parameter}, the resulting algorithm reading as:
\begin{align}
        u_{0}^{n+2/3} =\begin{cases}
          u^0 &\mbox{ if } n=0,\\
          u^{(n-1)+2/3} &\mbox{ if } n\geq 1.
          \label{eq.maevd.SQP.damp.0}
        \end{cases}
\end{align}
For $k\geq 0, u_k^{n+2/3}\rightarrow u_{k+1}^{n+2/3}$ as follows:
\begin{align}
&u_{k+1/2}^{n+2/3}=u^{n+1/3}+\nonumber\\
&\left[\frac{1-\left\|u_k^{n+2/3}\right\|_{L^{d+1}(\Omega)}^{d+1}+ (d+1){\displaystyle\int_{\Omega}} \left( u_k^{n+2/3}-u^{n+1/3}\right) u_k^{n+2/3} \left|u_k^{n+2/3}\right|^{d-1} d\x}{(d+1)\left\|u_k^{n+2/3}\right\|_{L^{2d}(\Omega)}^{2d}}\right] u_k^{n+2/3} \left| u_k^{n+2/3}\right|^{d-1}.
\label{eq.maevd.SQP.damp.1}
\end{align}
\begin{align}
 u_{k+1}^{n+2/3}=u_k^{n+2/3}+\beta_k\left( u_{k+1/2}^{n+2/3}-u_k^{n+2/3}\right) \ \left(=(1-\beta_k)u_k^{n+2/3}+\beta_k u_{k+1/2}^{n+2/3}\right),
 \label{eq.maevd.SQP.damp.2}
\end{align}
with $0<\delta\leq\beta_k\leq 1$. Remarks \ref{remark.spliterror}, \ref{remark.lambda} and \ref{remark.initial} still apply to algorithm (\ref{eq.maevd.split.0})-(\ref{eq.maevd.split.3}) with (\ref{eq.maev.initial}) replaced by
$$
u_0=\frac{\psi_0}{\|\psi_0\|_{L^{d+1}(\Omega)}} \mbox{ and } \p_0=\D^2 u_0.
$$
We will add that if the discrete analogue of algorithm (\ref{eq.maevd.split.0})-(\ref{eq.maevd.split.3}) is properly initialized, numerical experiments suggest that the discrete analogue of Step (\ref{eq.maevd.split.3}) is useless.

\subsection{Operator-splitting solution of problem (\ref{eq.mabg.ivp})}
The operator-splitting algorithm we employ to solve problem (\ref{eq.maevd.ivp}) reads as
\begin{align}
  (u^0,\p^0)=(u_0,\p_0).
  \label{eq.mabg.split.0}
\end{align}
For $n\geq0, (u^n,\p^n)\rightarrow u^{n+1/3}\rightarrow (u^{n+2/3},\p^{n+1})\rightarrow u^{n+1}$ as follows:\\
\underline{\emph{First Step}}: Solve the following \emph{(well-posed) linear elliptic problem}
\begin{align}
  \begin{cases}
    u^{n+1/3}-\tau \nabla\cdot \left[ \left(\varepsilon\I+\cof(\p^n)\right)\nabla u^{n+1/3}\right] =u^n \mbox{ in } \Omega,\\
    u^{n+1/3}=0 \mbox{ on }\partial\Omega.
  \end{cases}
  \label{eq.mabg.split.1}
\end{align}
\underline{\emph{Second Step}}:
\begin{subequations}
  \begin{align}
    &\p^{n+1}(\x)=P_+\left[ e^{-\gamma\tau} \p^n(\x)+\left(1-e^{-\gamma\tau}\right)\D^2 u^{n+1/3}(\x)\right], \mbox{ a.e. } \x\in\Omega, \label{eq.mabg.split.2a}\\
    &\begin{cases}
      u^{n+2/3}-u^{n+1/3}=-\tau d\lambda^{n+1} e^{-u^{n+2/3}},\\
      {\displaystyle \int_{\Omega}} \left(e^{-u^{n+2/3}}-1\right)d\x=C\ (>0)\ \left(\Leftrightarrow u^{n+2/3}\in \Sigma_{C}=\left\{ \phi \mbox{ measureble}, {\displaystyle \int_{\Omega}}\left(e^{-\phi}-1\right)d\x=C\right\}\right).
    \end{cases}
    \label{eq.mabg.split.2b}
  \end{align}
  \label{eq.mabg.split.2}
\end{subequations}
\underline{\emph{Third Step}}:
\begin{align}
  u^{n+1}=\inf (0,u^{n+2/3}).
  \label{eq.mabg.split.3}
\end{align}
Above,
\begin{itemize}
  \item It follows from (\ref{eq.mabg.split.1}) that $\left\| u^{n+1/3}\right\|_{L^2(\Omega)}\leq \left\| u^n\right\|_{L^2(\Omega)}, n>0$.
  \item Since we are looking for a \emph{ground state solution}, we will consider (\ref{eq.mabg.split.2b}) as an optimality system for the following minimization problem
      \begin{align}
        u^{n+2/3}\in \arg\min_{v\in\Sigma_C} \left( \frac{1}{2} \int_{\Omega} |v|^2d\x-\int_{\Omega} u^{n+1/3} vd\x\right).
        \label{eq.mabg.split.2.solve}
      \end{align}
      Unlike its finite element analogues, problem (\ref{eq.mabg.split.2.solve}) has no solution in general. However, as done in Section \ref{sec.timediscretize.maevd} for the solution of problem (\ref{eq.maev.split.2b.solve}), we will take advantage of the simpler formalism of the continuous problem to describe an algorithm of the SQP type for the formal solution of problem (\ref{eq.mabg.split.2.solve}). This SQP algorithm reads as
      \begin{align}
        u_{0}^{n+2/3} =\begin{cases}
          u^0 &\mbox{ if } n=0,\\
          u^{(n-1)+2/3} &\mbox{ if } n\geq 1.
          \label{eq.mabg.SQP.0}
        \end{cases}
      \end{align}
      For $k\geq 0, u_k^{n+2/3}\rightarrow u_{k+1}^{n+2/3}$ as follows:
      \begin{align}
        &u_{k+1}^{n+2/3}=u^{n+1/3}+\frac{\displaystyle\int_{\Omega} e^{-u_k^{n+2/3}}\left(1-u^{n+1/3}+u_k^{n+2/3}\right)d\x-(C+|\Omega|)}{\displaystyle\int_{\Omega} e^{-2u_k^{n+2/3}}d\x} e^{-u_k^{n+2/3}}.
        \label{eq.mabg.SQP.1}
      \end{align}
\end{itemize}
\begin{rmrk}
 As observed in Section \ref{sec.timediscretize.maev}, Remark \ref{remark.lambda}, there is no need to compute $\lambda^{n+1}$ in order to obtain the nonlinear eigenvalue $\lambda$. It follows from (\ref{eq.mabg.split.2b}) that $\lambda^{n+1}=\frac{\displaystyle\int_{\Omega}\left(u^{n+1/3}-u^{n+2/3}\right) d\x}{d\tau(C+|\Omega|)}$, the ratio of two small numbers. Following Remark \ref{remark.lambda}, we suggest approximating $\lambda$ by the following generalized \emph{Rayleigh quotient}
 \begin{align}
 \lambda_{\tau}=-\frac{\displaystyle\int_{\Omega} (\varepsilon\I+\cof(\p_{\tau}))\nabla u_{\tau} \cdot \nabla u_{\tau}d\x}{d\displaystyle\int_{\Omega} u_{\tau} e^{-u_{\tau}}d\x},
 \end{align}
 where $(u_{\tau},\p_{\tau})$ is the limit of the sequence $(u^{n},\p^n)_{n\geq 0}$.
\end{rmrk}
To conclude this section, let us mention that, as in Section \ref{sec.timediscretize.maevd}, if one initializes properly algorithm (\ref{eq.mabg.split.0})-(\ref{eq.mabg.split.3}) (an issue we will address in Section \ref{sec.numerical}), Step (\ref{eq.mabg.split.3}) is useless.
\section{Finite Element Approximation of the Nonlinear Eigenvalue Problems}
\label{sec.spacediscretize}
\subsection{Generalities}
\begin{figure}[t]
  \centering
  \includegraphics[width=0.2\textwidth]{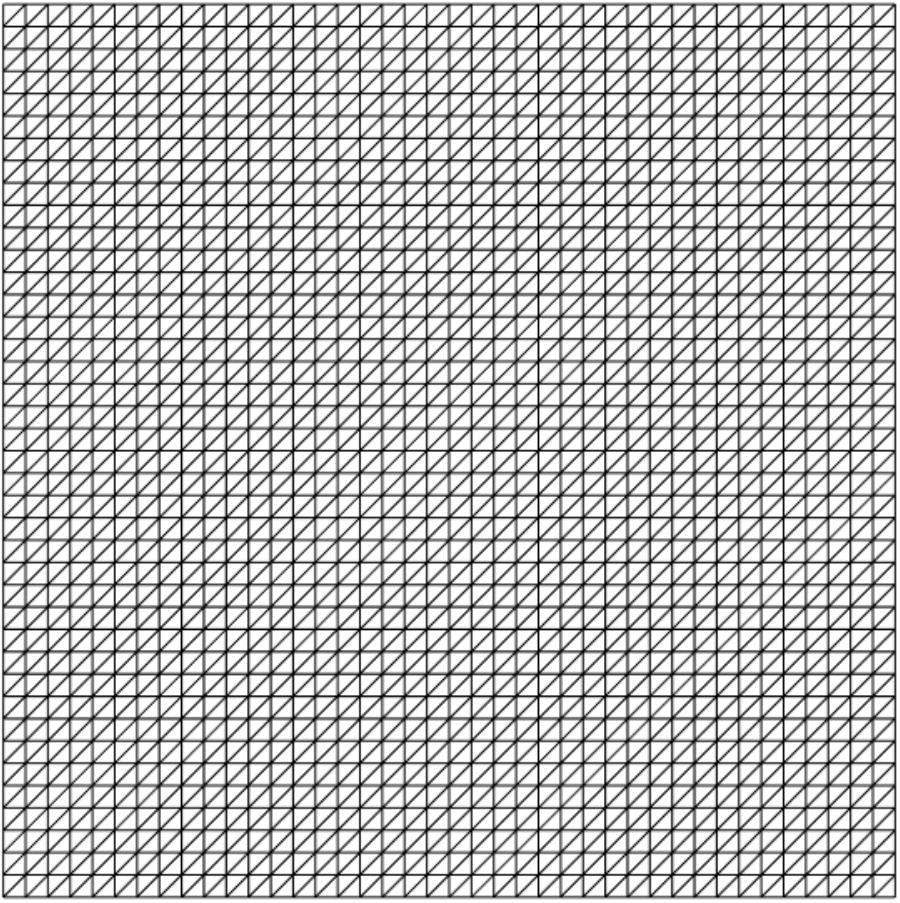}\hspace{1cm}
  \includegraphics[width=0.2\textwidth]{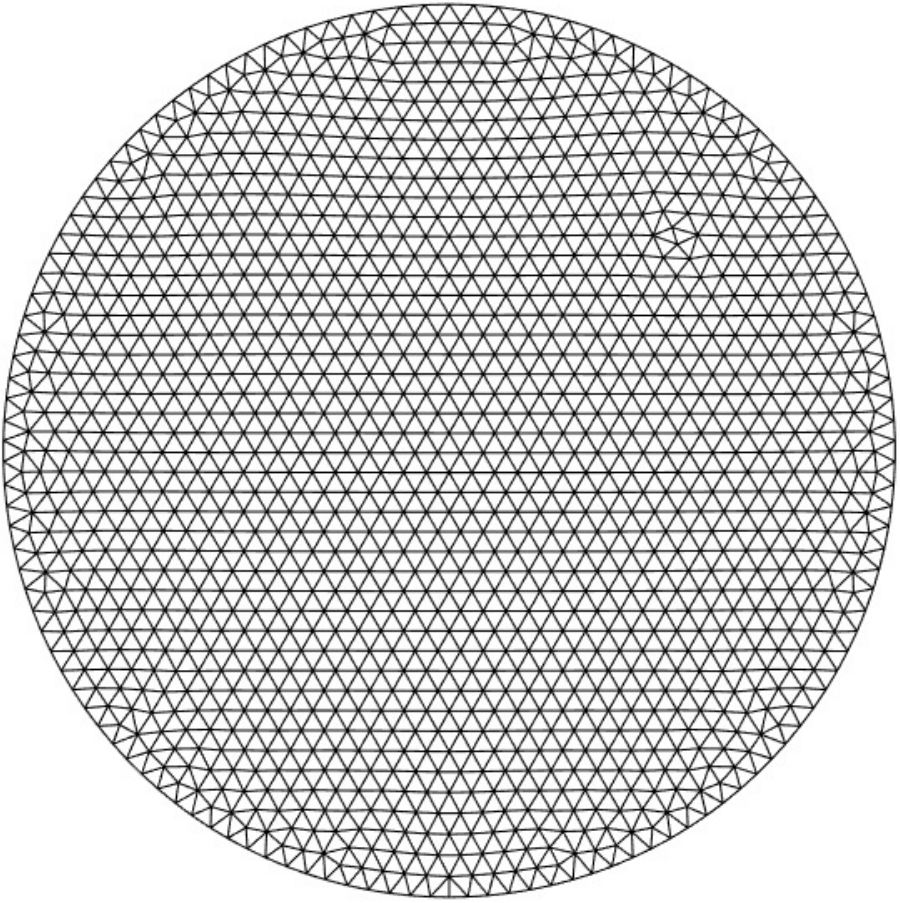}
  \caption{Finite element triangulations of a square and of a disk}\label{fig.mesh}
\end{figure}
Assuming that $\Omega$ is a \emph{bounded convex polygonal domain} of $\RR^2$ (or has been approximated by a family of such domains) we introduce a family $(\mT_h)_h$ of triangulations of $\Omega$ like those in Figure \ref{fig.mesh} ($h$ is, typically, the \emph{length of the largest edge(s)} of $\mT_h$). Next, we approximate the functional spaces $H^1(\Omega)$ and $H_0^1(\Omega)$ by
\begin{align}
  V_h=\left\{ \phi\in C^0(\bar{\Omega}), \phi|_{T} \in \bP_1, \forall T\in \mT_h\right\},
\end{align}
and
\begin{align}
  V_{0h}=\left\{ \phi\in V_h, \phi|_{\partial\Omega}=0\right\}\ (=V_h\cap H_0^1(\Omega)),
\end{align}
respectively, $\bP_1$ being the space of the polynomial functions of two variables of degree$\leq1$. Let us denote by $\Sigma_h$ (resp., $\Sigma_{0h}$) the set of the \emph{vertices} of $\mT_h$ (resp., the set $\Sigma_h\backslash\Sigma_h \cap \partial\Omega$). We have then
$$
\dim V_h =\mathrm{Card}\ \Sigma_h\ (=N_h) \mbox{ and } \dim V_{0h}=\mathrm{Card}\ \Sigma_{0h}\ (=N_{0h}).
$$
We assume that the vertices of $\mT_h$ have been numbered so that $\Sigma_{0h}=\{Q_k\}_{k=1}^{N_{0h}}$. For $k=1,...,N_h$, we define $\omega_k$ as the union of those triangles of $\mT_h$ which have $Q_k$ as a common vertex. We denote by $|\omega_k|$ the measure (area) of $\omega_k$.

Problem (\ref{eq.mabg.ivp}) being the more complicated of the initial value problems we considered in Section \ref{sec.intro}, \ref{sec.formulation} and \ref{sec.timediscretize}, we will focus on the finite element implementation of algorithm (\ref{eq.mabg.split.0})-(\ref{eq.mabg.split.3}). The resulting algorithm is easy to modify in order to handle problems (\ref{eq.maevd.ivp}) and (\ref{eq.maev.ivp}).

\subsection{On the finite element approximation of problem (\ref{eq.mabg.split.1})}
First, we equip space $V_h$ with the inner product $\{\phi,\theta\}\rightarrow (\phi,\theta)_h: V_h\times V_h\rightarrow \RR$ defined by
\begin{align}
  (\phi,\theta)_h=\frac{1}{3} \sum_{k=1}^{N_h} |\omega_k|\phi(Q_k)\theta(Q_k), \forall \phi,\theta\in V_h.
\end{align}
Next, we approximate the elliptic problem (\ref{eq.mabg.split.1}) by
\begin{align}
  \begin{cases}
    u_h^{n+1/3}\in V_{0h},\\
    \left( u_h^{n+1/3},v\right)_h+ \tau \displaystyle\int_{\Omega} (\varepsilon\I+\cof(\p_h^n))\nabla u_h^{n+1/3}\cdot \nabla v d\x=(u_h^n,v)_h,\\
    \forall v\in V_{0h}.
  \end{cases}
  \label{eq.mabf.dis.1}
\end{align}
If the matrix $\p_h^n(Q_k)$ is \emph{positive semi-definite}, $\forall k=1,...,N_h$, then problem (\ref{eq.mabf.dis.1}) has a \emph{unique} solution. Indeed, problem (\ref{eq.mabf.dis.1}) is equivalent to a linear system associated with a sparse symmetric positive definite matrix. In \cite{glowinski2019finite}, one addresses the direct and iterative solution of discrete elliptic problems such as (\ref{eq.mabf.dis.1}).

\subsection{Approximating $\p^{n+1}$}
Proceeding as in \cite{glowinski2019finite}, we define $\p_h^{n+1}$, the approximation of $\p^{n+1}$ in relation (\ref{eq.mabg.split.2a}), by:
\begin{align}
  \begin{cases}
    \p^{n+1}_h\in (V_h)^{2\times 2},\\
    \p_h^{n+1}(Q_k)=P_+\left[ e^{-\gamma\tau}\p_h^{n}(Q_k)+\left(1-e^{-\gamma\tau}\right) \D_h^2 u_h^{n+1/3}(Q_k)\right],\\
    \forall k=1,...,N_h,
  \end{cases}
  \label{eq.mabf.dis.2}
\end{align}
where in (\ref{eq.mabf.dis.2}) the approximate Hessian $\D_h^2 u_h^{n+1/3}$ is obtained as follows:\\
For $1\leq i,j\leq 2$, solve
\begin{align}
\begin{cases}
  \pi_{ijh}^{n+1/3}\in V_{0h},\\
  c\displaystyle\sum_{T\in\mT_h} |T|\displaystyle\int_T \nabla \pi_{ijh}^{n+1/3}\cdot \nabla \phi d\x+ \left(\pi_{ijh}^{n+1/3},\phi\right)_h= -\frac{1}{2} \displaystyle\int_{\Omega}\left[ \frac{\partial u_h^{n+1/3}}{\partial x_i} \frac{\partial \phi}{\partial x_j}+ \frac{\partial u_h^{n+1/3}}{\partial x_j} \frac{\partial \phi}{\partial x_i} \right]d\x,\\
  \forall \phi\in V_{0h},
\end{cases}
\label{eq.mabf.dis.3}
\end{align}
and
\begin{align}
  \begin{cases}
    D_{ijh}^2 u_h^{n+1/3}\in V_h,\\
    c\displaystyle\sum_{T\in\mT_h} |T|\displaystyle\int_T \nabla D_{ijh}^2 u_h^{n+1/3}\cdot \nabla \phi d\x+ \left(D_{ijh}^2 u_h^{n+1/3},\phi\right)_h= \left(\pi_{ijh}^{n+1/3},\phi\right)_h,\\
    \forall \phi\in V_h,
  \end{cases}
  \label{eq.mabf.dis.4}
\end{align}
with $c\approx 1$ in (\ref{eq.mabf.dis.3}) and (\ref{eq.mabf.dis.4}), and set
\begin{align}
  \D_h^2u_h^{n+1/3}=\left( D_{ijh}^2u_h^{n+1/3}\right)_{1\leq i,j\leq 2}.
\end{align}
The numerical experiments reported in reference \cite{glowinski2019finite} shows that the \emph{double regularization method} based on (\ref{eq.mabf.dis.3}) and (\ref{eq.mabf.dis.4}) is well-suited to Monge-Amp\`{e}re equations with non-smooth solutions.

\subsection{Approximating $u^{n+2/3}$}
Let us denote by $u_h^{n+2/3}$ the approximation of $u^{n+2/3}$. To compute $u_h^{n+2/3}$ we employ the following discrete analogue of the SQP algorithm (\ref{eq.mabg.SQP.0}), (\ref{eq.mabg.SQP.1}):
\begin{align}
  u_{0h}^{n+1/3}=\begin{cases}
    u_h^0, &\mbox{ if } n=0,\\
    u_h^{(n-1)+2/3} &\mbox{ if } n\geq 1.
  \end{cases}
  \label{eq.mabf.dis.5}
\end{align}
For $q\geq0, u_{qh}^{n+2/3}\rightarrow u_{(q+1)h}^{n+2/3}$ as follows
\begin{align}
  \begin{cases}
    u_{(q+1)h}^{n+2/3}(Q_k)=u_h^{n+1/3}(Q_k) + \frac{\displaystyle\sum_{l=1}^{N_{0h}}|\omega_l|e^{-u_{qh}^{n+2/3}(Q_l)}\left( 1-u_h^{n+1/3}+u_{qh}^{n+2/3}\right)(Q_l)-3(C+|\Omega|)}{\displaystyle\sum_{l=1}^{N_{0h}} |\omega_l|e^{-2u_{qh}^{n+2/3}(Q_l)}} e^{-u_{qh}^{n+2/3}}(Q_k),\\
    \forall k=1,...,N_{0h}.
  \end{cases}
  \label{eq.mabf.dis.6}
\end{align}
Since $\tau$ is small, an obvious variant of algorithm (\ref{eq.mabf.dis.5}), (\ref{eq.mabf.dis.6}) is the one obtained by performing only one iteration of the above algorithm.
\section{Numerical Experiments}
\label{sec.numerical}
\subsection{Generalities}
\label{sec.numerical.general}
In order to validate the methodology discussed in the above sections we performed various numerical experiments concerning the three nonlinear eigenvalue problems we introduced in Section \ref{sec.intro} (namely, problem (\ref{eq.maevd}), (\ref{eq.maev}) and (\ref{eq.mabg})). These experiments include the cases where $\Omega$ is the unit disk of $\RR^2$, the simplicity of the geometry allowing a very accurate computation of the solutions, and therefore meaningful comparisons. Proceeding by increasing order of complexity, we will start our result presentation with problem (\ref{eq.maev}) and will conclude with problem (\ref{eq.mabg}), by far the most challenging.

\begin{table}[t]
  \centering
  \begin{tabular}{|c|c|c|c|c|c|c|c|c|}
    \hline
    $h$& Iteration \# & $\|u_h^{n+1}-u_h^n\|_{0h}$ & $L^2$-error & rate &$L^{\infty}$-error &rate&$\lambda_h$& $\min u_h$ \\
    \hline
    1/10 & 248 & 9.66$\times 10^{-10}$ & $6.18\times10^{-2}$ & & $5.73\times10^{-2}$& &3.64 & -0.9639\\
    \hline
    1/20 & 774 & 9.93$\times 10^{-10}$ &$4.74\times10^{-2}$ & 0.38& $4.24\times10^{-2}$ & 0.43& 4.52 & -0.9857\\
    \hline
    1/40 & 2520 & 9.97$\times 10^{-10}$ & $2.86\times10^{-2}$ & 0.73& $2.76\times10^{-2}$& 0.62 &5.10 & -1.0033\\
    \hline
    1/80 & 7427 & 9.99$\times 10^{-10}$ & $1.54\times10^{-2}$ & 0.92& $1.58\times10^{-2}$ & 0.80&5.40 & -1.0134\\
    \hline
    1/160 & 31345 & 9.99$\times 10^{-11}$ & $7.58\times10^{-3}$ & 1.02& $8.57\times10^{-3}$ & 0.88&5.56 & -1.0189\\
    \hline
  \end{tabular}
  \caption{Problem (\ref{eq.maev.num}) with $\Omega=\left\{ \x=\{x_1,x_2\},x_1^2+x_2^2<1\right\}$. Variations with $h$ of the number of iterations necessary to achieve convergence (2nd column), of the $L^2$ and $L^{\infty}$ approximation errors and of the associated convergence rates (columns 4, 5, 6 and 7), of the computed eigenvalue (8th column) and of the minimal value of $u_h$ over $\Omega$ (that is $u_h(\mathbf{0})$) (9th column).}
  \label{tab.maev.disk}
\end{table}

\begin{figure}[t]
  (a)
  \includegraphics[width=0.35\textwidth]{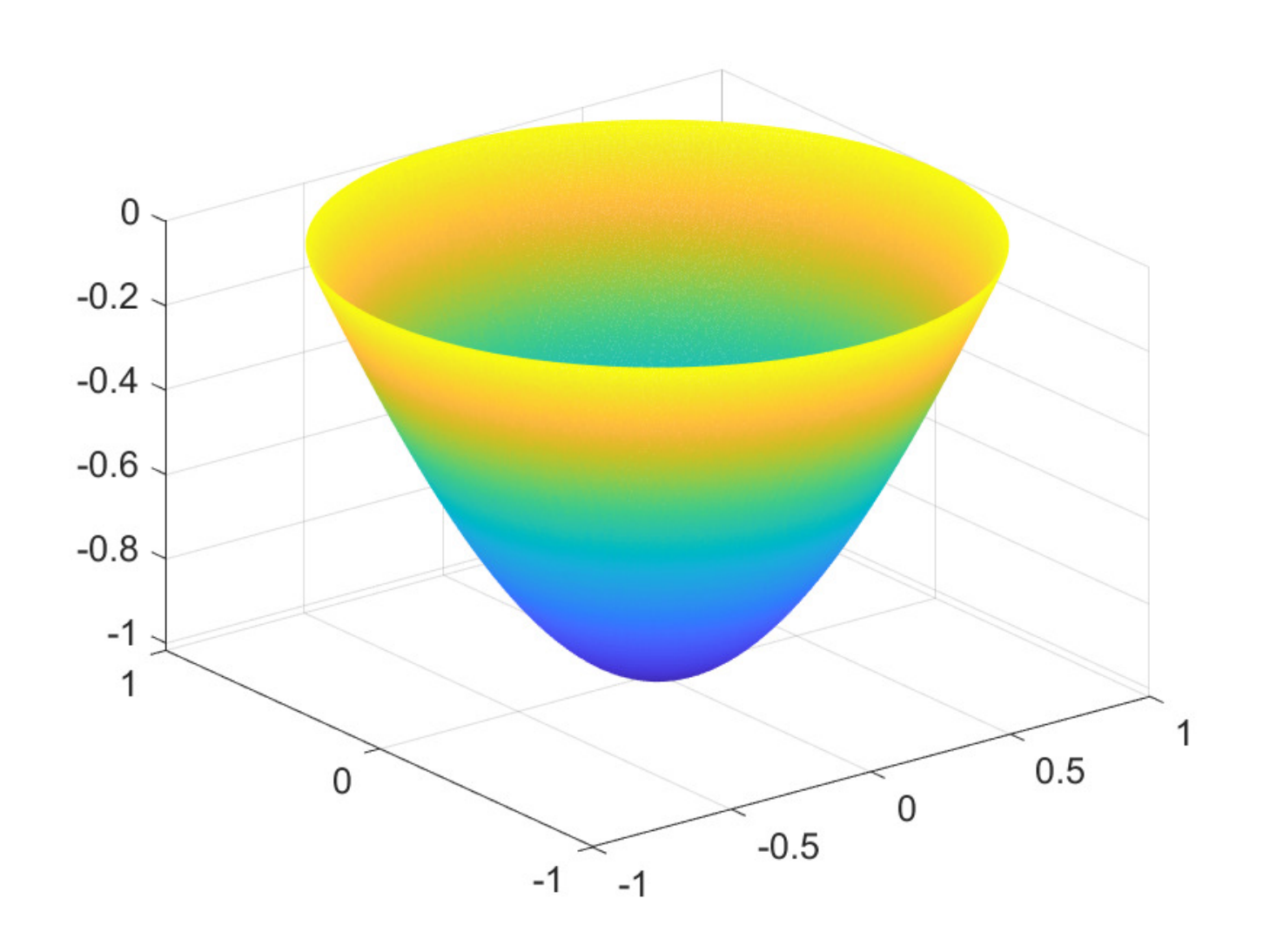}
  (b)
  \includegraphics[width=0.35\textwidth]{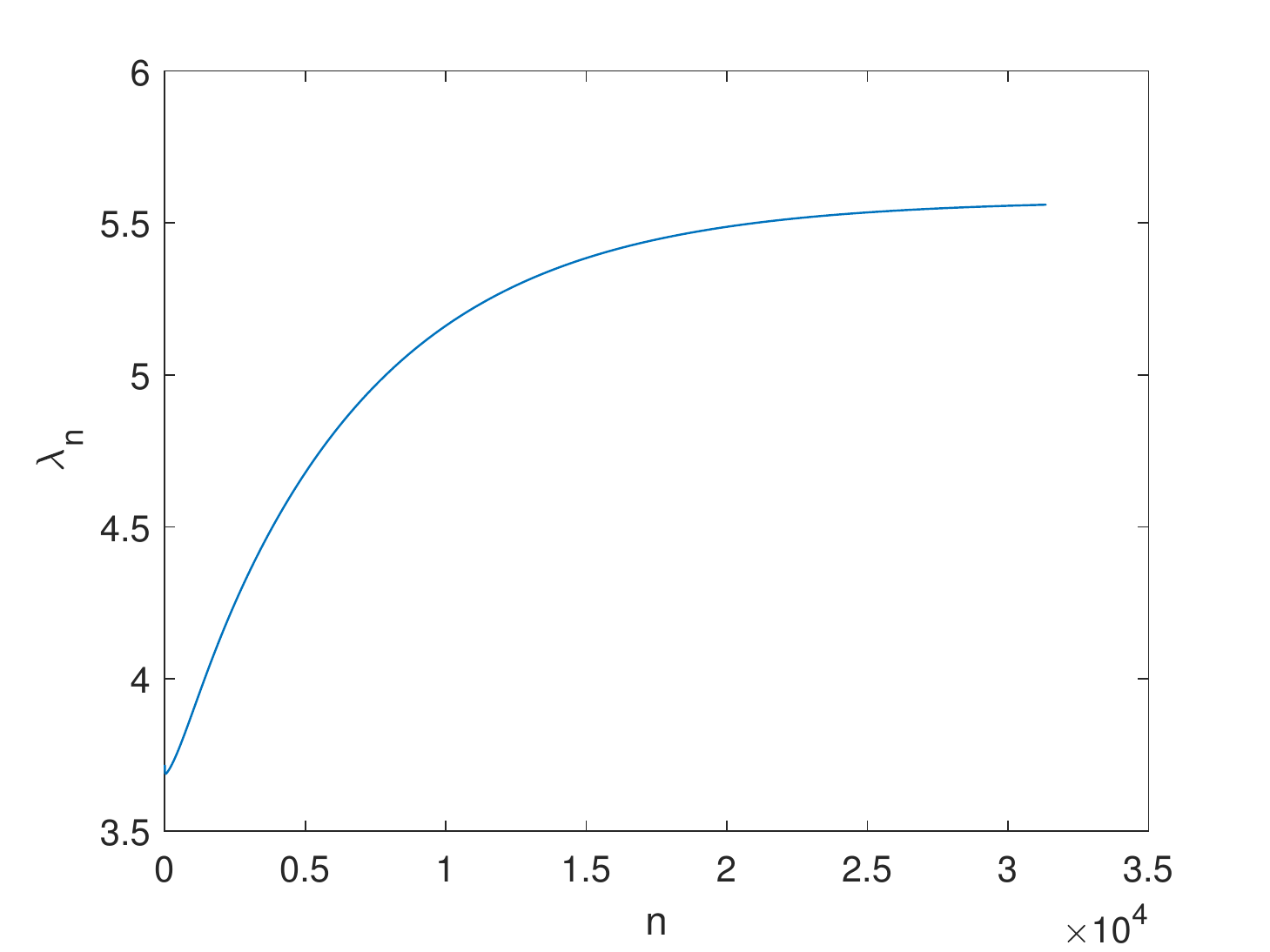}\\
  (c)
  \includegraphics[width=0.35\textwidth]{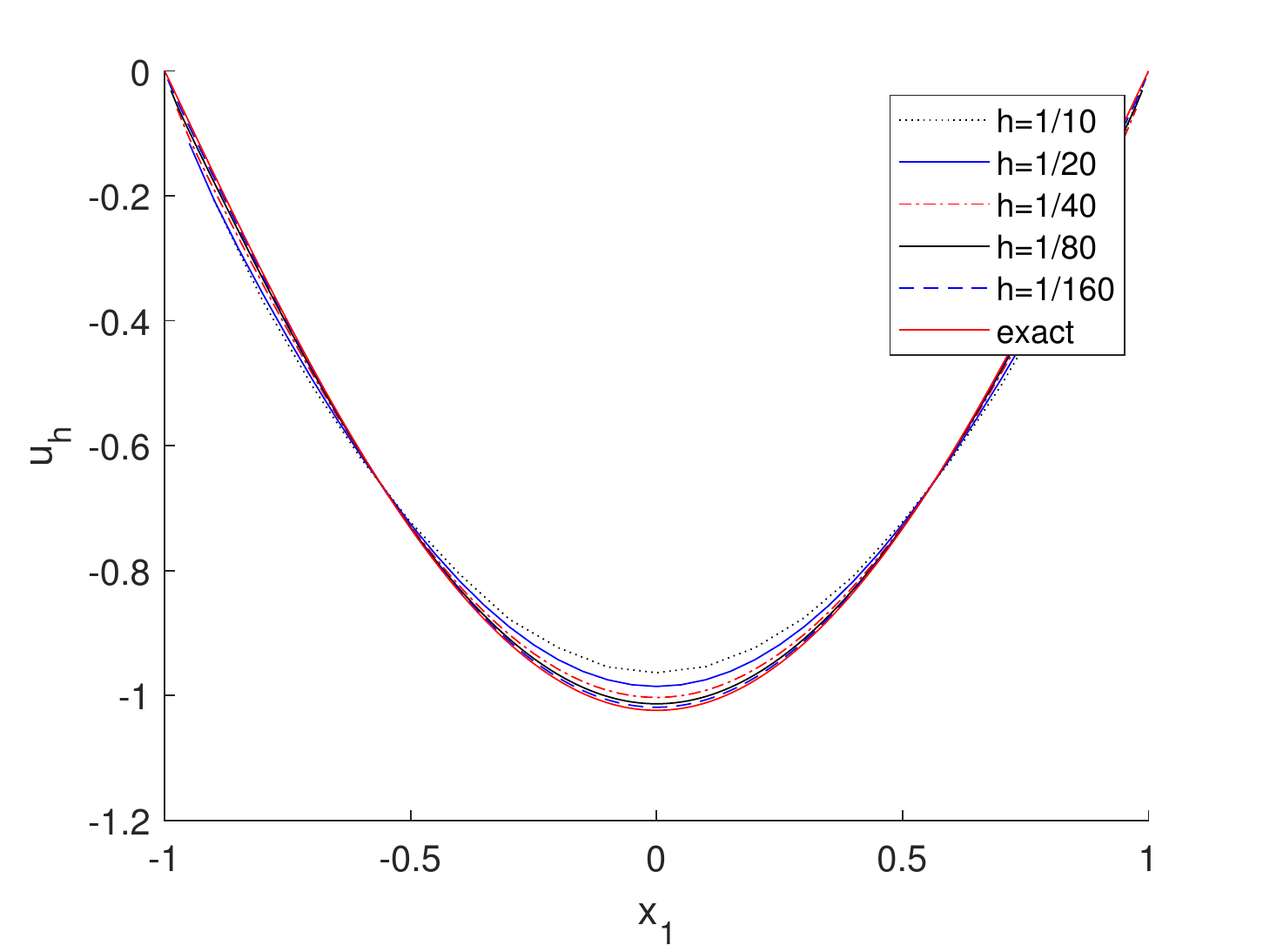}
  (d)
  \includegraphics[width=0.35\textwidth]{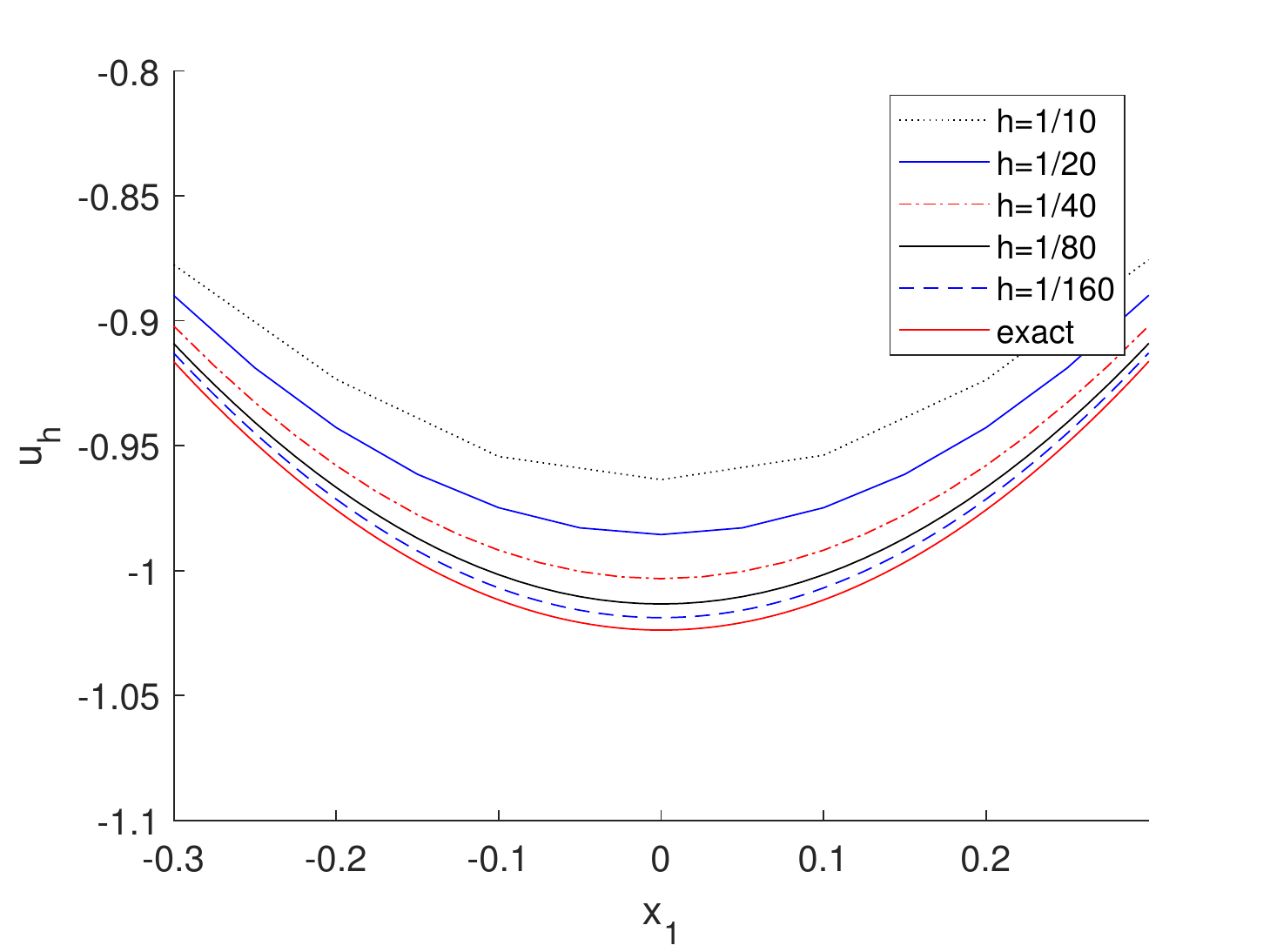}\\
  \caption{Problem (\ref{eq.maev.num}) with $\Omega=\left\{ \x=\{x_1,x_2\},x_1^2+x_2^2<1\right\}$. (a) Graph of the approximate solution $u_h$ computed with $h=1/160$. (b) Variation with $n$ of the discrete Rayleigh quotient approximating $\lambda_n$ for $h=1/160$. (c) Graphs of the restrictions of $u_h$ and $u$ to the line $x_2=0$ for $h=1/10,1/20,1/40,1/80$ and $1/160$. (d) Zoom version of (c).}
  \label{fig.maev.disk}
\end{figure}

\begin{figure}[t]
  \centering
  \includegraphics[width=0.3\textwidth]{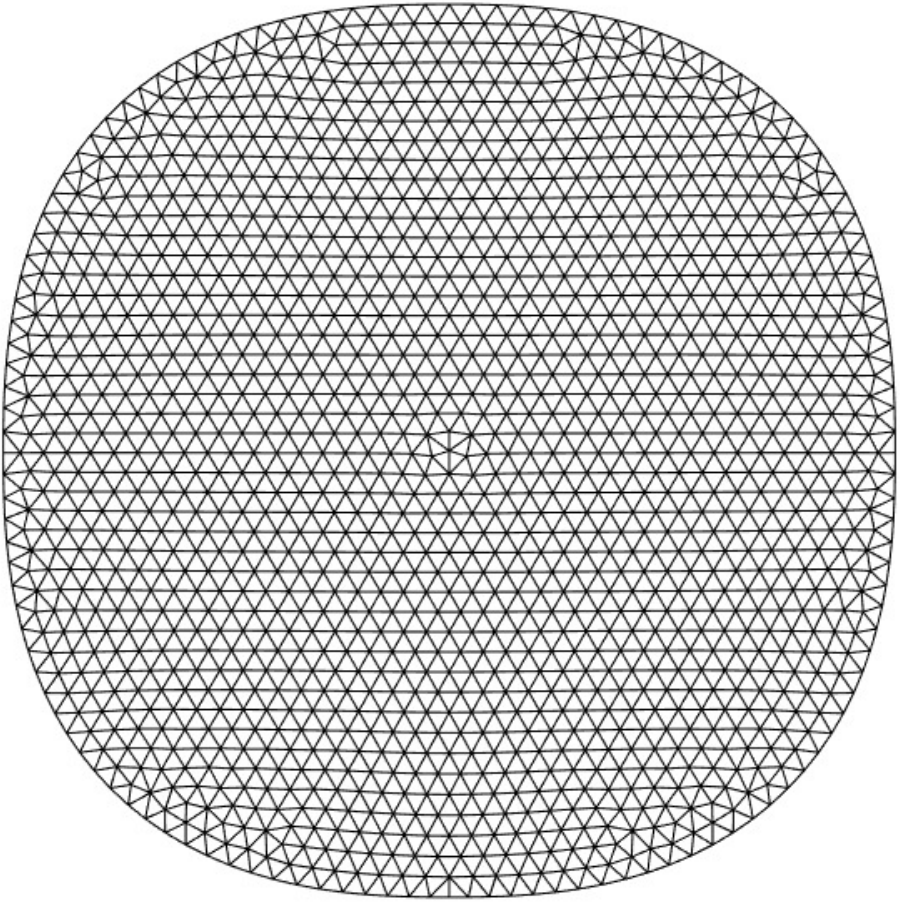}
  \caption{Isotropic finite element triangulation of the convex domain $\Omega$ defined by (\ref{eq.square2o5}) ($h=1/20$).}\label{fig.square2o5}
\end{figure}

\begin{table}[t]
  \centering
  \begin{tabular}{|c|c|c|c|c|}
    \hline
    $h$& Iteration \# & $\|u^{n+1}_h-u_h^n\|_{0h}$ & $\lambda_h$& $\min u_h$\\
    \hline
    1/10 & 233 & 9.67$\times 10^{-10}$ & 3.07 & -0.9286\\
    \hline
    1/20 & 807 & 9.96$\times 10^{-10}$ & 3.78 & -0.9506\\
    \hline
    1/40 & 2648 & 9.97$\times 10^{-10}$ & 4.22 & -0.9660\\
    \hline
    1/80 & 7981 & 9.99$\times 10^{-10}$ & 4.45 & -0.9754\\
    \hline
  \end{tabular}
  \caption{Problem (\ref{eq.maev.num}) with $\Omega$ defined by (\ref{eq.square2o5}). Variations with $h$ of the number of iterations necessary to achieve convergence (2nd column), of the computed eigenvalue (4th column) and of the minimal value of $u_h$ over $\Omega$ (that is $u_h(\mathbf{0})$) (5th column).}
  \label{tab.maev.square2o5}
\end{table}

\begin{figure}[t]
  (a)
  \includegraphics[width=0.35\textwidth]{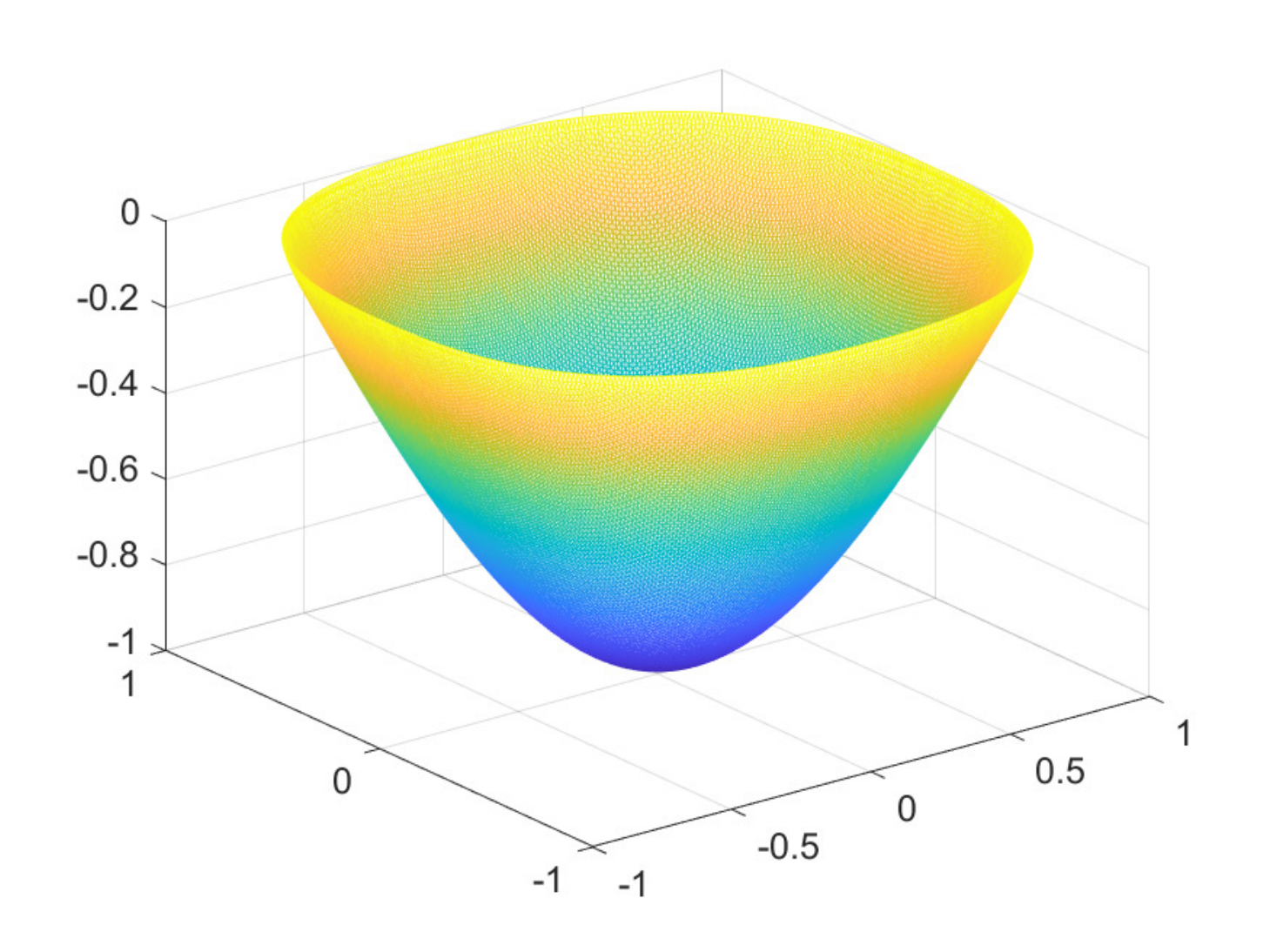}
  (b)
  \includegraphics[width=0.35\textwidth]{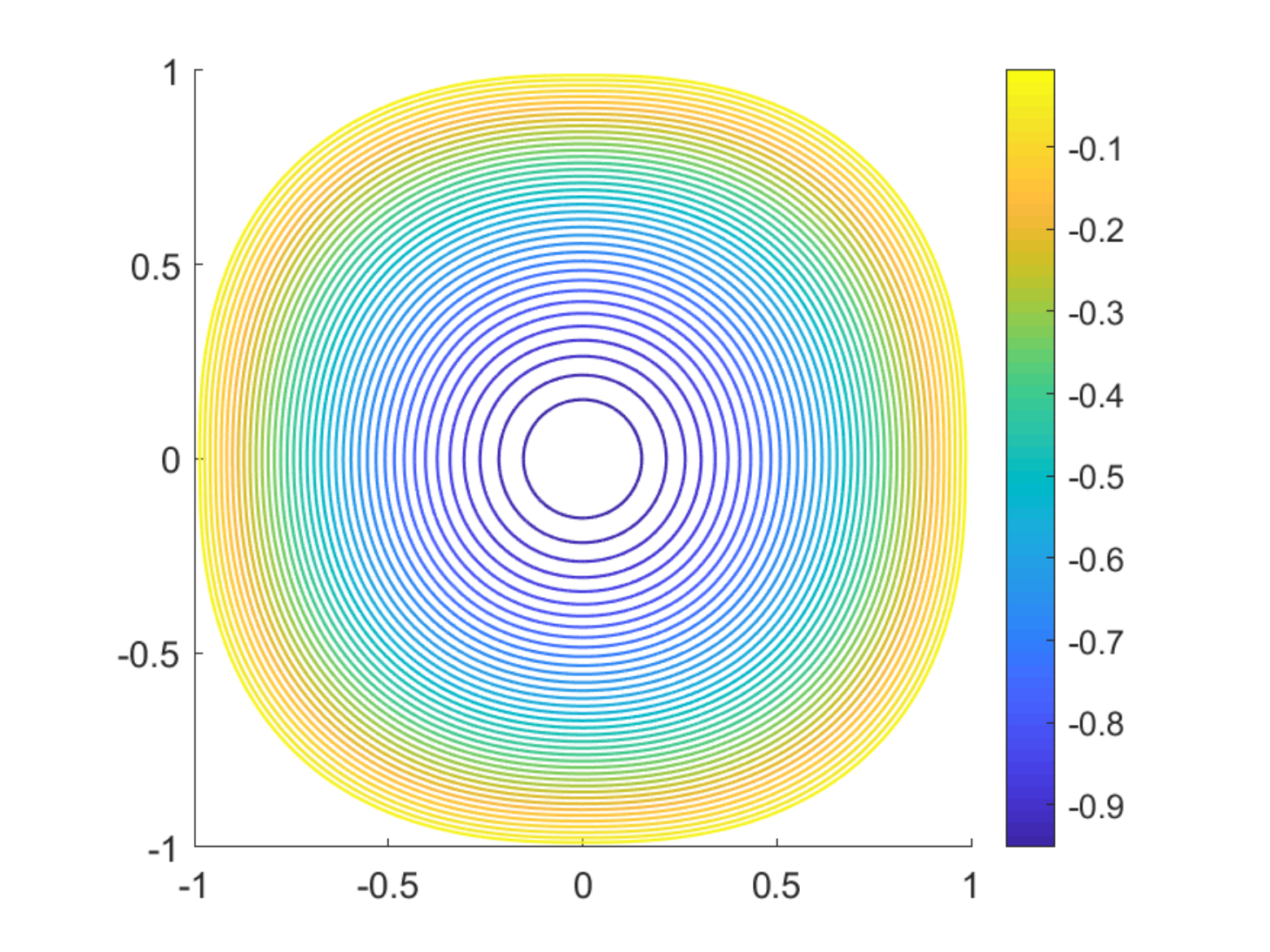}\\
  (c)
  \includegraphics[width=0.35\textwidth]{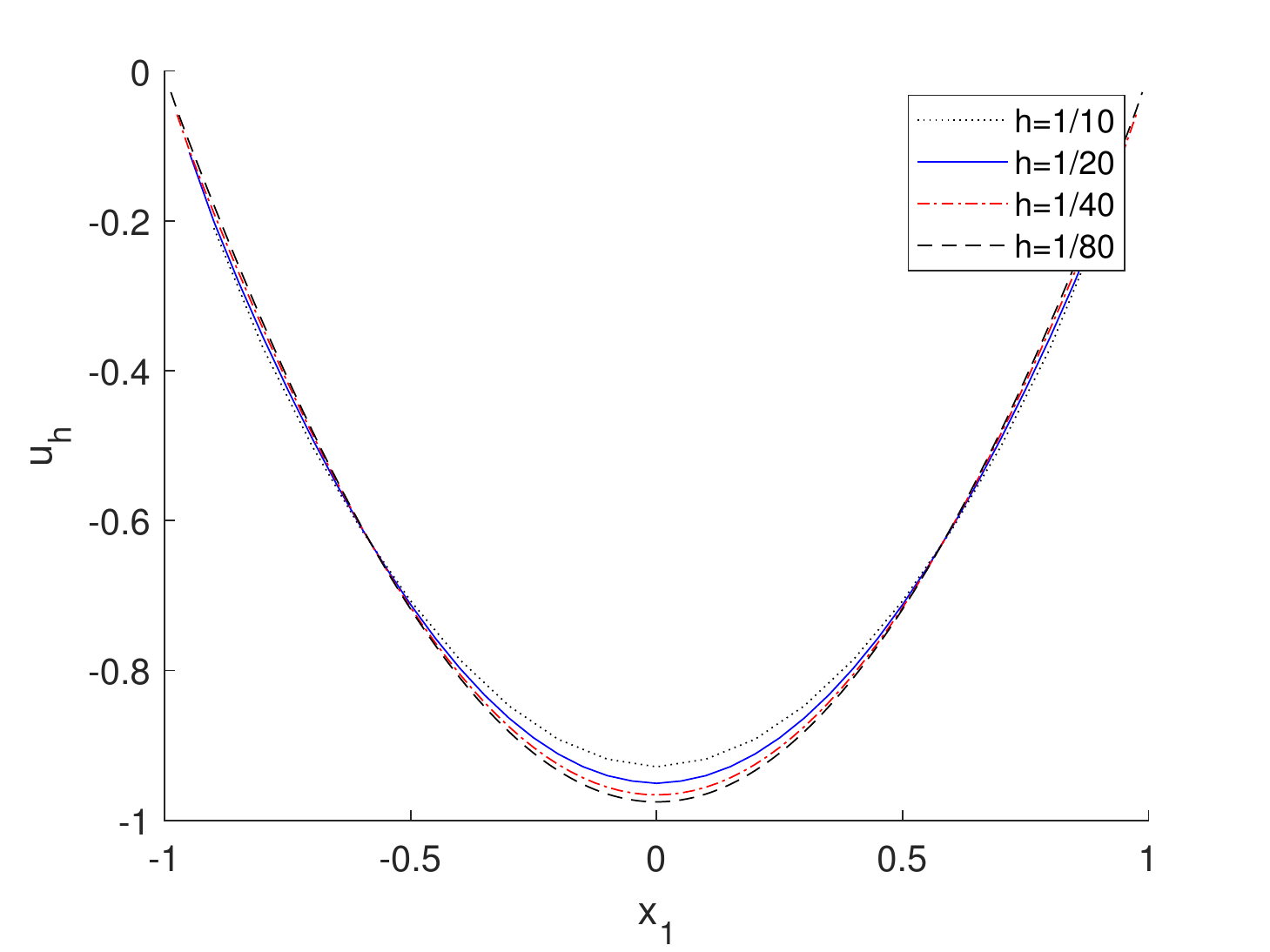}
  \caption{Problem (\ref{eq.maev.num}) with $\Omega$ defined by (\ref{eq.square2o5}): (a) Graph of the approximate solution $u_h$ computed with $h=1/80$. (b) Contour of the approximated solution $u_h$ for $h=1/80$. (c) Graphs of the restrictions of $u_h$ and $u$ to the line $x_2=0$ for $h=1/10,1/20,1/40$ and $1/80$.}
  \label{fig.maev.square2o5}
\end{figure}
\subsection{On the solution of problem (\ref{eq.maev})}
From Section \ref{sec.intro} and \ref{sec.formulation}, the particular problem (\ref{eq.maev}) we consider is to find the \emph{ground state solution} of
\begin{align}
  \begin{cases}
    u\leq0, \lambda>0,\\
    \det\D^2u=-\lambda u \mbox{ in } \Omega,\\
    u=0 \mbox{ on } \partial\Omega,\\
    \displaystyle\int_{\Omega}|u|^2d\x=1.
  \end{cases}
  \label{eq.maev.num}
\end{align}
As announced in Section \ref{sec.numerical.general}, we investigate first the situation where in (\ref{eq.maev.num}), $\Omega$ is the \emph{unit disk} of $\RR^2$ centered at $\mathbf{0}$, that is $\Omega=\left\{\x=\{x_1,x_2\},x_1^2+x_2^2<1\right\}$. It makes sense to assume that in that special case, problem (\ref{eq.maev.num}) has a smooth radial solution, this solution verifying:
\begin{align}
  \begin{cases}
    u\leq0, \lambda\geq0,\\
    u'u''=-\lambda r u \mbox{ in }(0,1),\\
    u'(0)=0, u(1)=0,\\
    2\pi\displaystyle\int_0^1 |u|^2rdr=1,
  \end{cases}
  \label{eq.maev.disk}
\end{align}
with $r=\sqrt{x_1^2+x_2^2}$. Using a \emph{shooting technique}, one can solve the two-point ODE problem (\ref{eq.maev.disk}) very accurately, obtaining a radial solution of (\ref{eq.maev.num}) verifying $u(\mathbf{0})=-1.0238...$ and $\lambda=5.7183...$.

The algorithms discussed in Section \ref{sec.timediscretize} and \ref{sec.spacediscretize} were implemented using \emph{isotropic unstructured finite element triangulations}, like the one on the right of Figure \ref{fig.mesh}, with $h$ ranging from $1/10$ to $1/160$. The results reported in Table \ref{tab.maev.disk} and Figure \ref{fig.maev.disk} have been obtained with $\varepsilon=\tau=h^2$, using $\|u_h^{n+1}-u_h^n\|_{0h}\leq 10^{-9}$ (resp., $\leq 10^{-10}$) if $h=1/10,1/20,1/40,1/80$ (resp., $h=1/160$) as stopping criterion (with $\|v\|_{0h}=\sqrt{\displaystyle\sum_{k=1}^{N_h} |\omega_k||v(Q_k)||^2}\ \left(=\sqrt{(v,v)_h}\right),\forall v\in V_h$).

The results reported in Table \ref{tab.maev.disk} suggest that the $L^2$-norm of  approximation error (in fact $\|u_h-u\|_{0h}$) is close to $O(h)$. Actually, it follows from the results reported in column 2, that the number of iterations necessary to achieve convergence varies like $h^{-\frac{\log 3}{\log 2}}$, approximately. On the other hand, the results from the 8th column show that $\lambda-\lambda_h=25h$, quite accurately. Additional results have been visualized on Figure \ref{fig.maev.disk}.

The next problem (\ref{eq.maev.num}) is the one associated with the bounded convex domain $\Omega$ defined by
\begin{align}
  \Omega=\left\{ \x=\{x_1,x_2\}, |x_1|^{2.5}+|x_2|^{2.5}<1\right\}.
  \label{eq.square2o5}
\end{align}
The algorithms discussed in Section \ref{sec.timediscretize} and \ref{sec.spacediscretize} were implemented using \emph{isotropic unstructured finite element triangulations}, like the one in Figure \ref{fig.square2o5}, with $h$ ranging from $1/10$ to $1/80$. The results reported in Table \ref{tab.maev.square2o5} and Figure \ref{fig.maev.square2o5} have been obtained with $\varepsilon=\tau=h^2$, using $\|u_h^{n+1}-u_h^n\|_{0h}\leq 10^{-9}$ as stopping criterion.

Using least-squares fitting leads us to
\begin{align}
\lambda_h\approx \lambda_0-ch
\end{align}
with $\lambda_0=4.61437...$ and $c=15.666...$. On Figure \ref{fig.maev.square2o5}, we have reported further information concerning the approximate solutions. The commonalities with the results obtained for the unit disk are striking but expected.

\begin{table}[t]
  \centering
  \begin{tabular}{|c|c|c|c|c|c|c|c|c|}
    \hline
    $h$& Iteration \# & $\|u^{n+1}_h-u_h^n\|_{0h}$ & $L^2$-error & rate & $L^{\infty}$-error & rate &$\lambda_h$&  $\min u_h$ \\
    \hline
    1/20 & 312 & 9.92$\times 10^{-7}$ & 4.20$\times 10^{-2}$ & & 4.02$\times 10^{-2}$  & & 6.13 &  -1.1405\\
    \hline
    1/40 & 976 & 9.99$\times 10^{-7}$ & 2.84$\times 10^{-2}$ & 0.56 & 2.81$\times 10^{-2}$  & 0.52 & 6.81 &   -1.1461\\
    \hline
    1/80& 1017 & 9.95$\times 10^{-7}$ & 1.43$\times 10^{-2}$ & 0.99 & 1.57$\times 10^{-2}$  & 0.84 & 7.12 & -1.1538\\
    \hline
    1/160 & 3315 & 9.99$\times 10^{-7}$ & 7.09$\times 10^{-3}$ & 1.01 & 8.74$\times 10^{-3}$  & 0.85 & 7.32 &  -1.1559\\
    \hline
  \end{tabular}
  \caption{Problem (\ref{eq.maevd.num}) with $\Omega=\left\{\x=\{x_1,x_2\},x_1^2+x_2^2<1\right\}$. Variations with $h$ of the number of iterations necessary to achieve convergence (2nd column), of the $L^2$ and $L^{\infty}$ approximation errors and of the associated convergence rates (column 4, 5, 6 and 7), of the computed eigenvalue (8th column) and of the minimal value of $u_h$ over $\Omega$ (that is $u_h(\mathbf{0})$) (9th column).}
  \label{tab.maevd.disk}
\end{table}

\begin{figure}[t]
  (a)
  \includegraphics[width=0.35\textwidth]{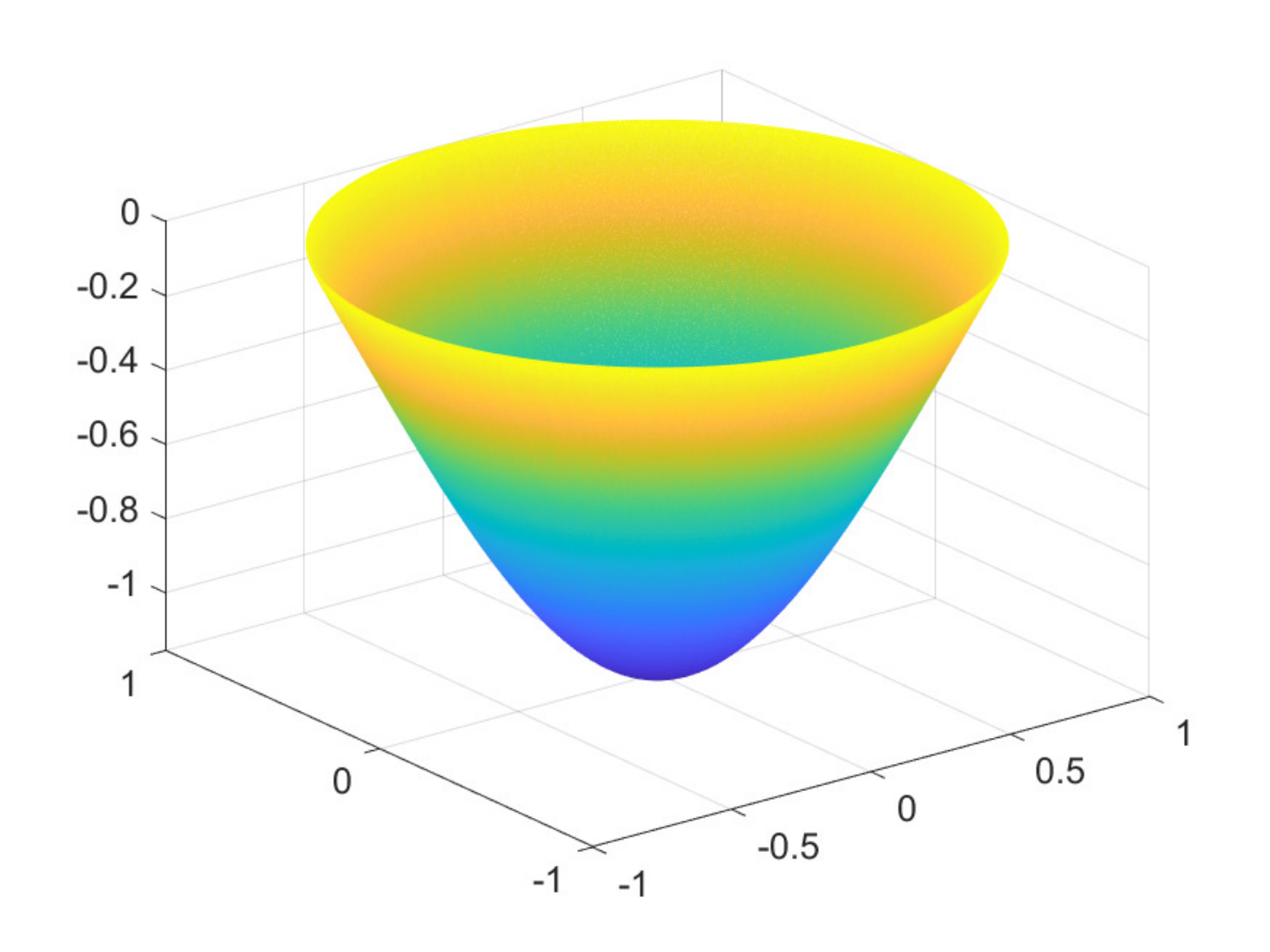}
  (b)
  \includegraphics[width=0.35\textwidth]{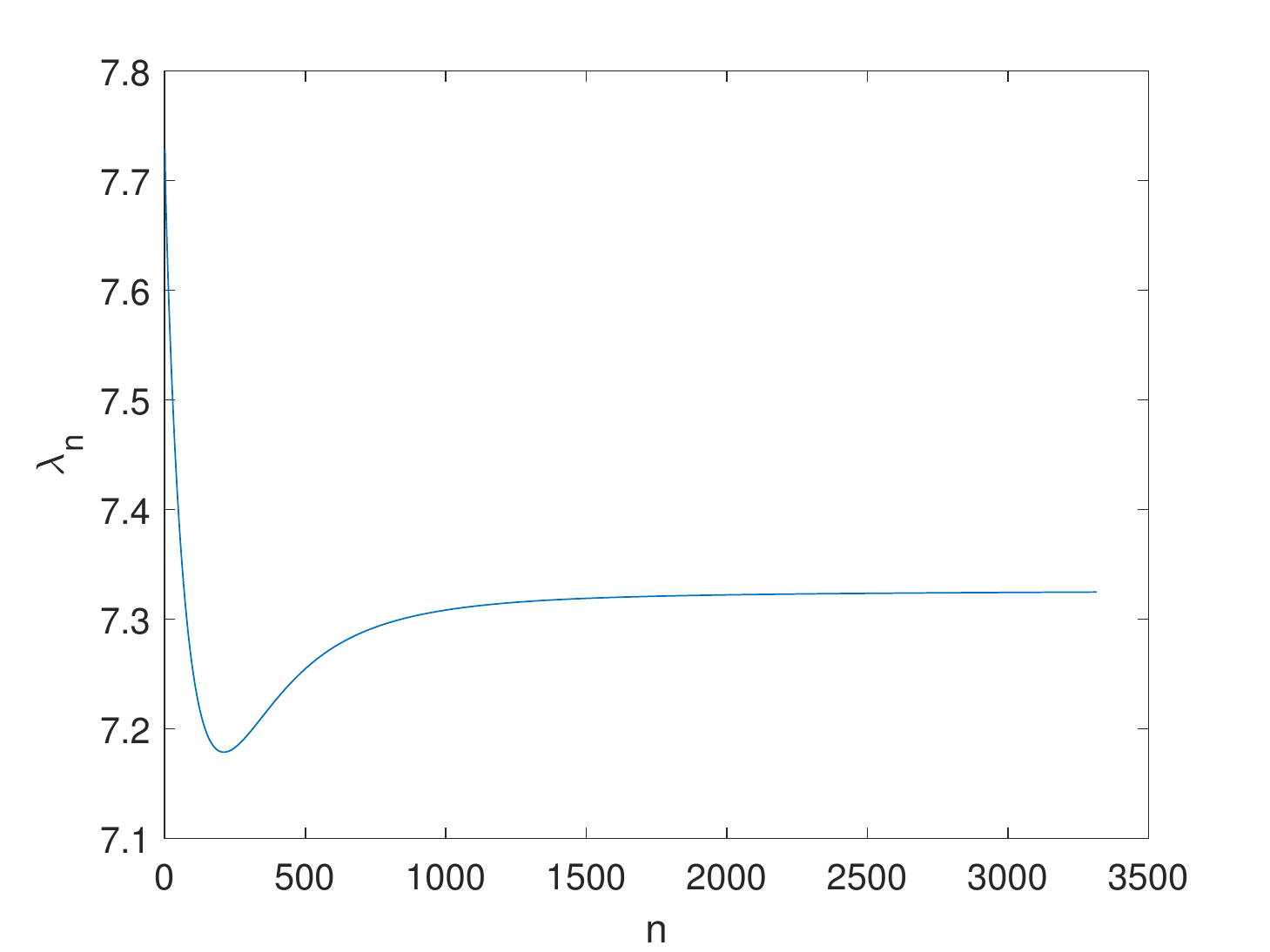}\\
  (c)
  \includegraphics[width=0.35\textwidth]{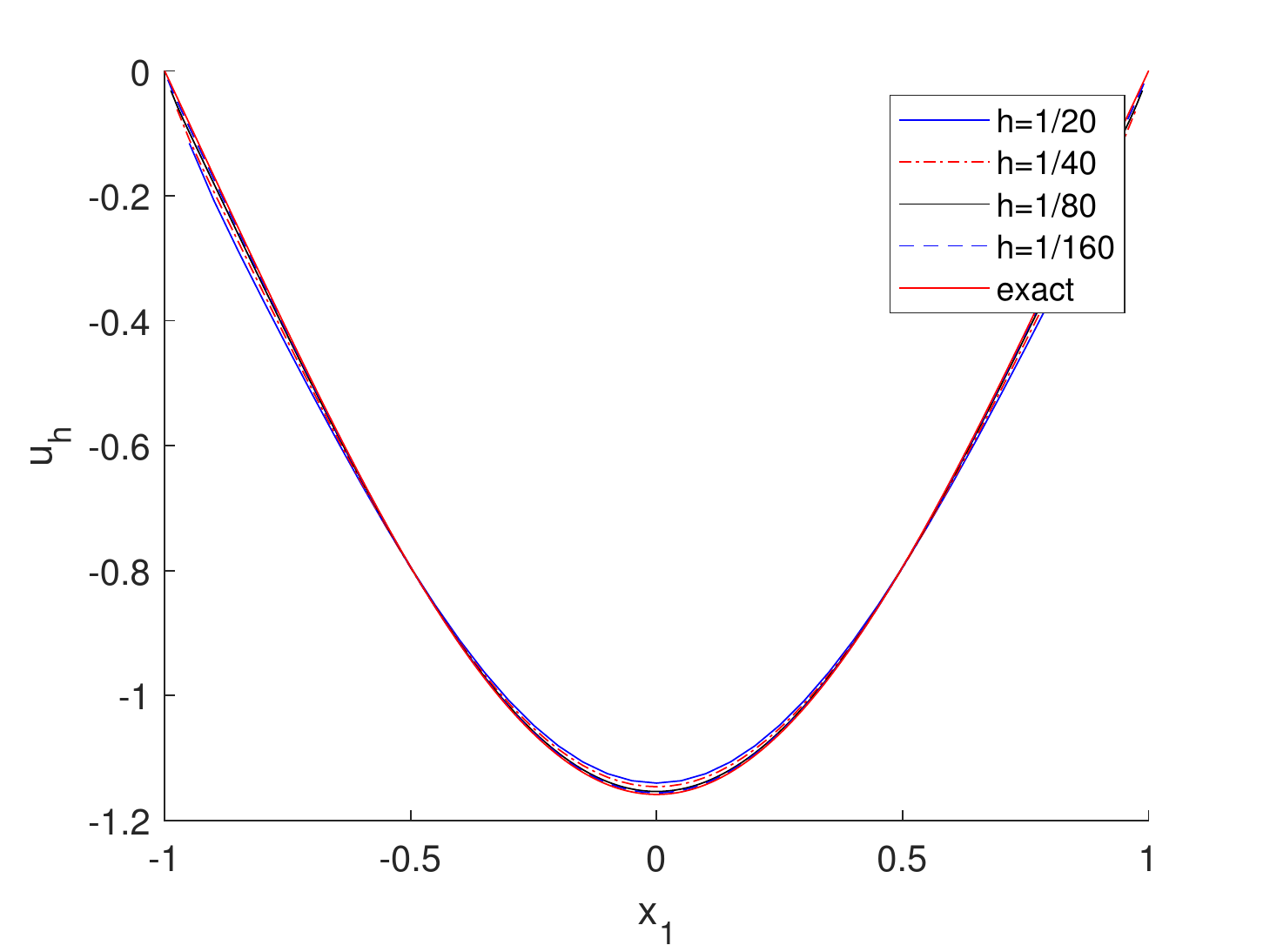}
  (d)
  \includegraphics[width=0.35\textwidth]{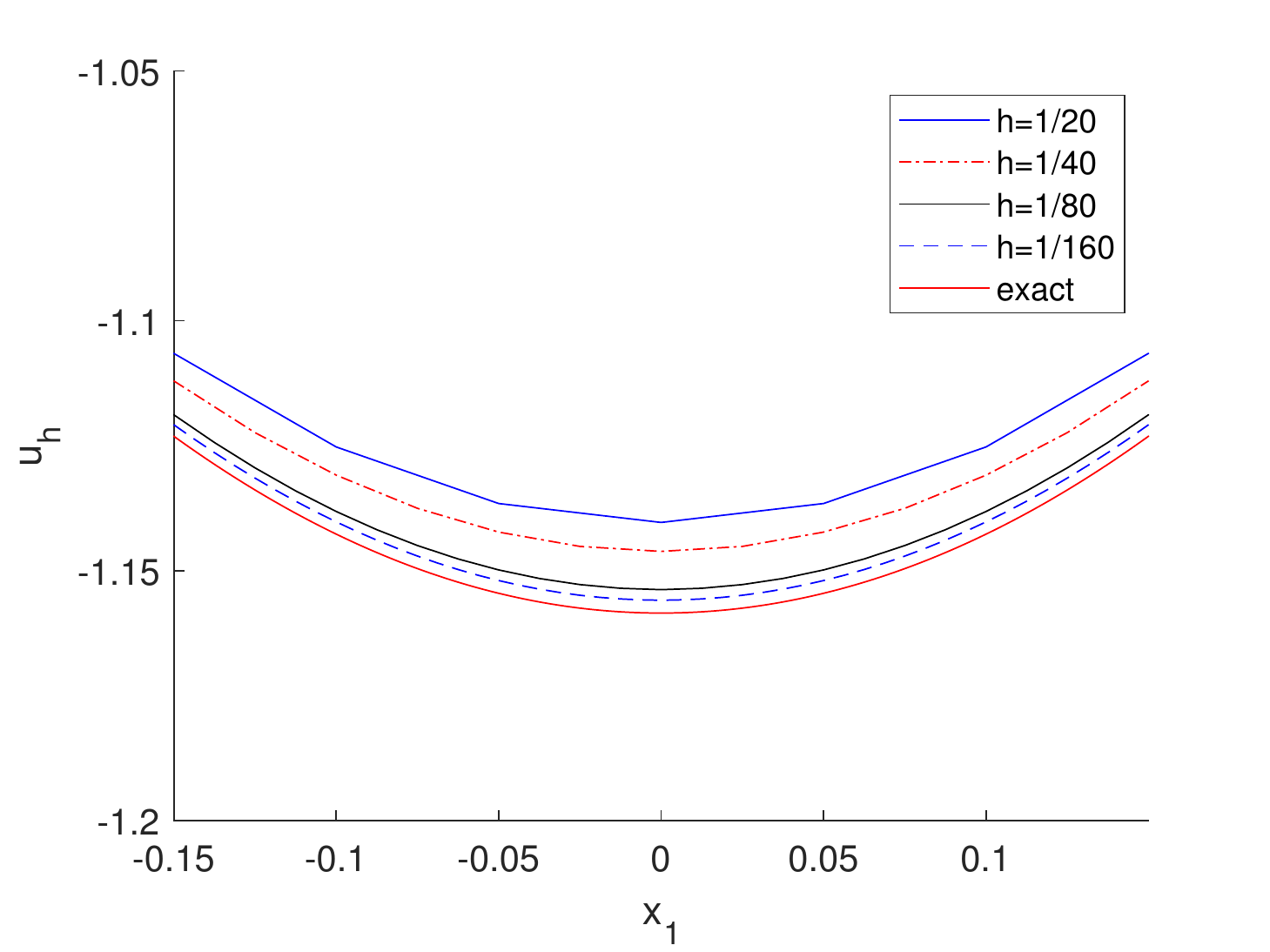}
  \caption{Problem (\ref{eq.maevd.num}) with $\Omega=\left\{ \x=\{x_1,x_2\},x_1^2+x_2^2<1\right\}$. (a) Graph of the approximate solution $u_h$ computed with $h=1/160$. (b) Variation with $n$ of the discrete Rayleigh quotient approximating $\lambda_n$ for $h=1/160$. (c) Graphs of the restrictions of $u_h$ and $u$ to the line $x_2=0$ for $h=1/20,1/40,1/80$ and $1/160$. (d) Zoom version of (c).}
  \label{fig.maevd.disk}
\end{figure}

\begin{table}[t]
  \centering
  \begin{tabular}{|c|c|c|c|c|}
    \hline
    $h$& Iteration \# & $\|u_h^{n+1}-u_h^n\|_{0h}$ & $\lambda_h$&  $\min u_h$\\
    \hline
    1/10 & 208 & 9.96$\times 10^{-9}$ & 4.35 & -1.1187\\
    \hline
    1/20 & 671 & 9.93$\times 10^{-9}$ & 5.30 & -1.1122\\
    \hline
    1/40 & 2164 & 9.98$\times 10^{-9}$ & 5.85 & -1.1170\\
    \hline
    1/80 & 6584 & 9.99$\times 10^{-9}$ & 6.13 & -1.1222\\
    \hline
  \end{tabular}
   \caption{Problem (\ref{eq.maevd.num}) with $\Omega$ defined by (\ref{eq.square2o5}). Variations with $h$ of the number of iterations necessary to achieve convergence (2nd column), of the computed eigenvalue (4th column) and of the minimal value of $u_h$ over $\Omega$ (that is $u_h(\mathbf{0})$) (5th column).}
  \label{tab.maevd.square2o5}
\end{table}

\begin{figure}[h]
  (a)
  \includegraphics[width=0.35\textwidth]{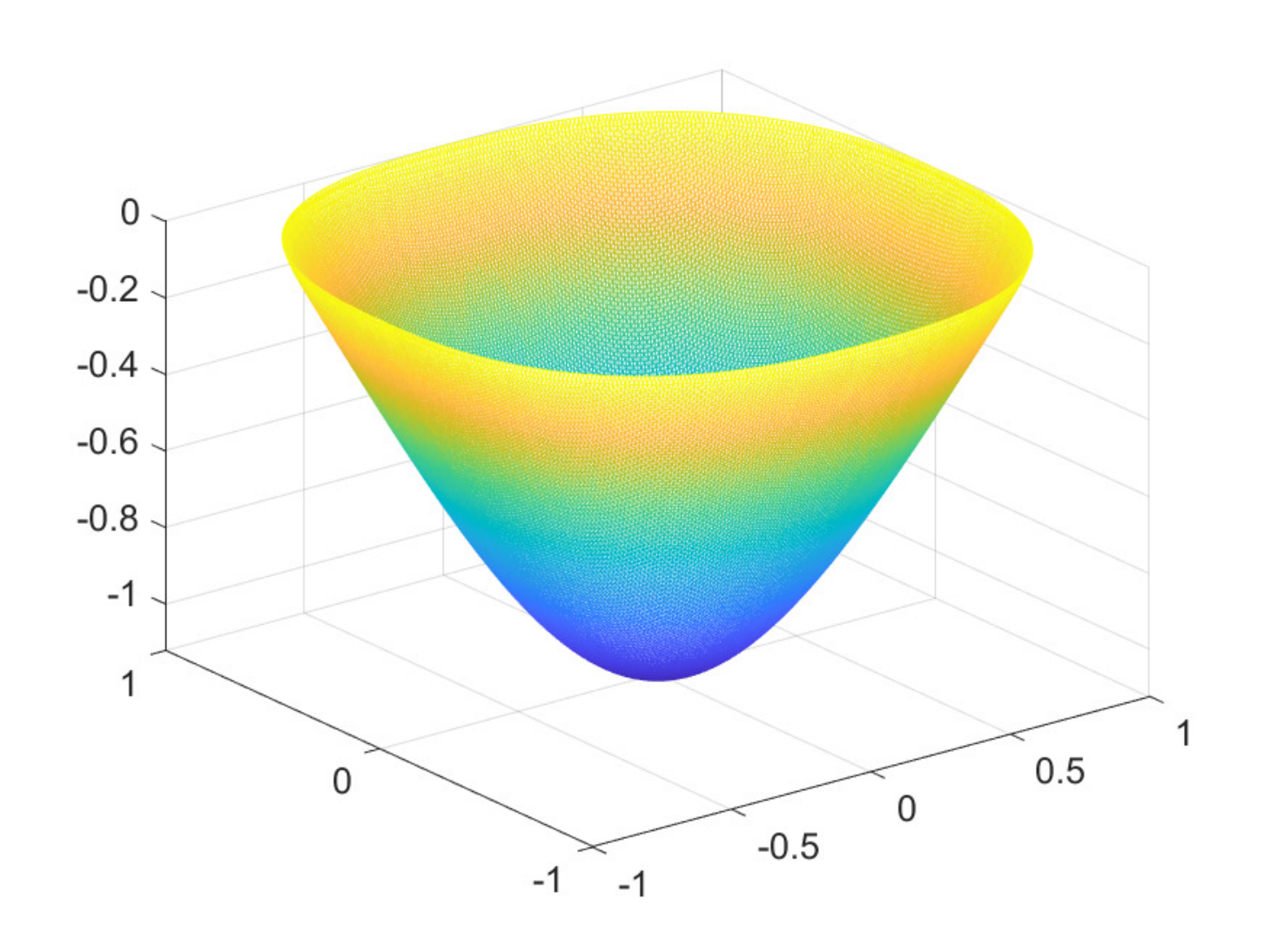}
  (b)
  \includegraphics[width=0.35\textwidth]{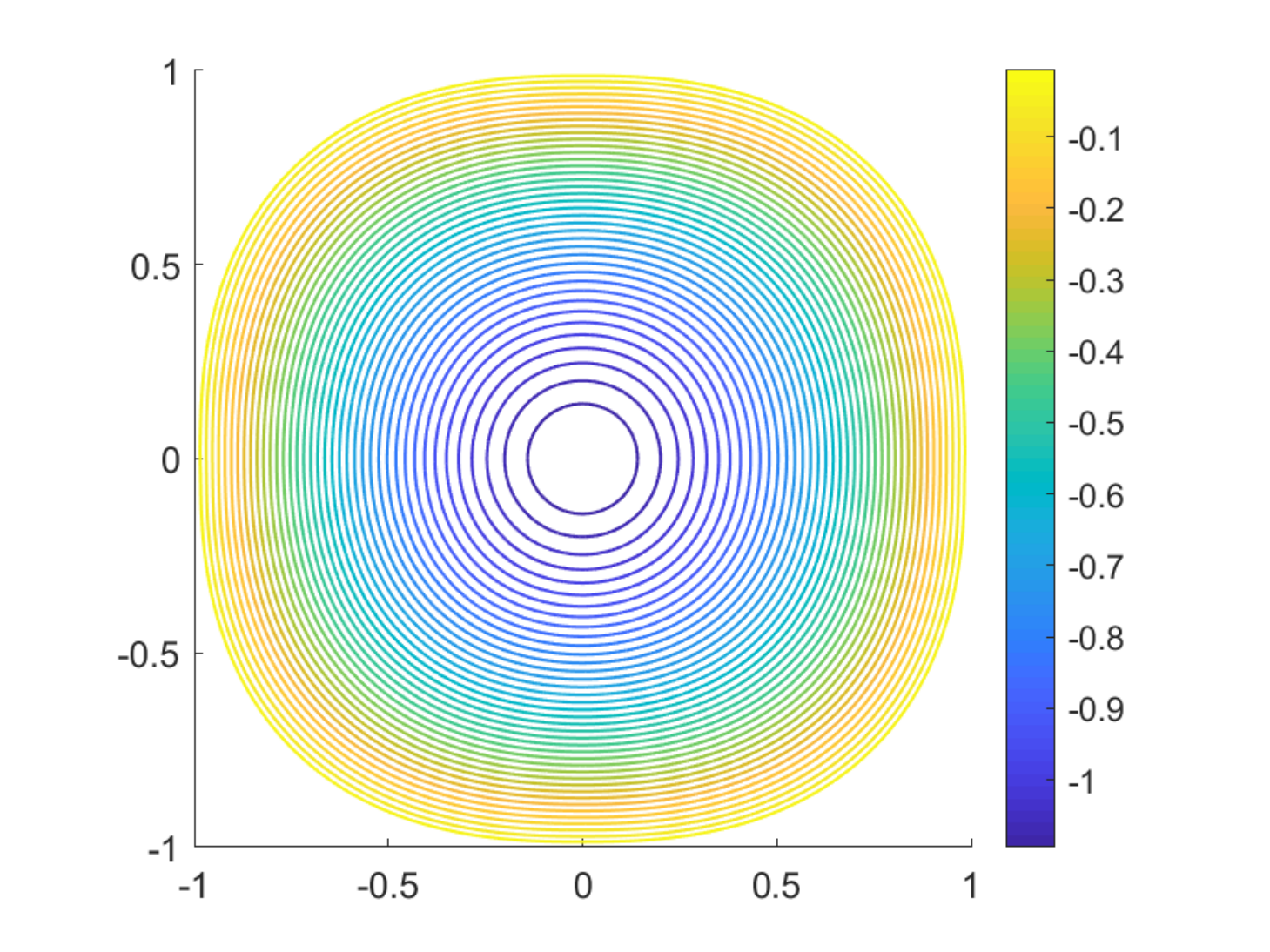}\\
  (c)
  \includegraphics[width=0.35\textwidth]{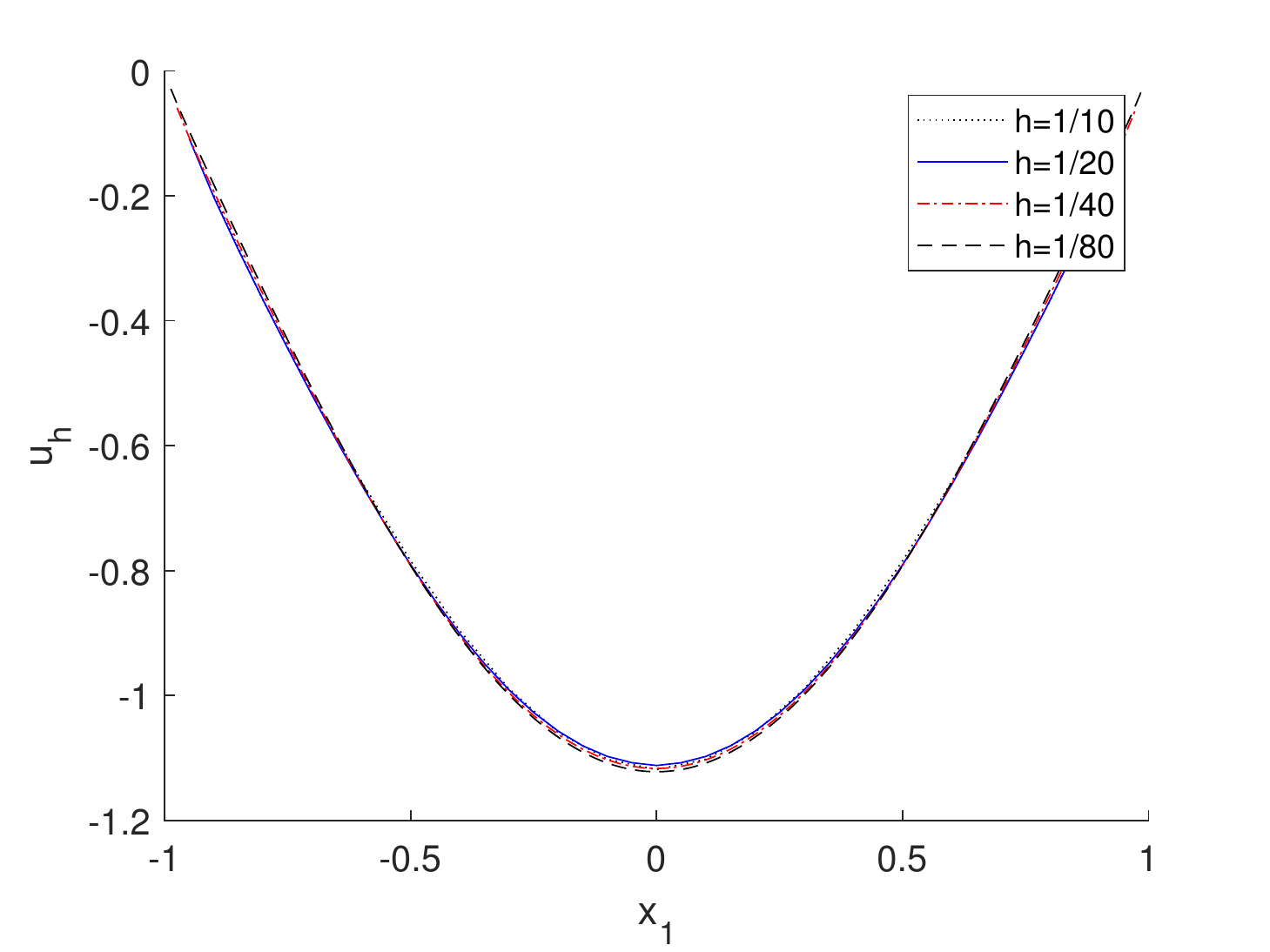}
  (d)
  \includegraphics[width=0.35\textwidth]{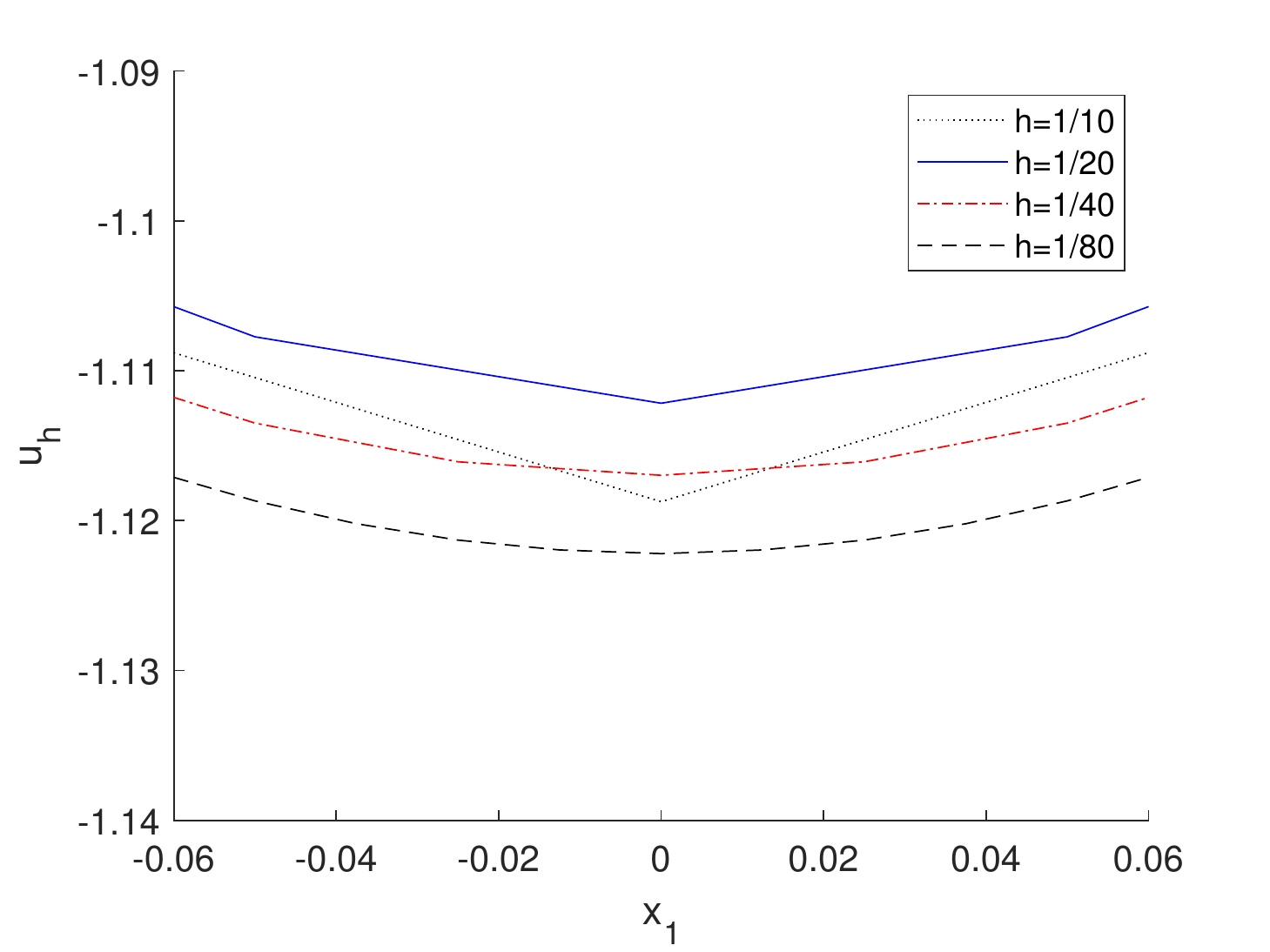}
  \caption{Problem (\ref{eq.maevd.num}) with $\Omega$ defined by (\ref{eq.square2o5}): (a) Graph of the approximate solution $u_h$ computed with $h=1/80$. (b) Contour of the approximated solution $u_h$ for $h=1/80$. (c) Graphs of the restrictions of $u_h$ and $u$ to the line $x_2=0$ for $h=1/10,1/20,1/40$ and $1/80$. (d) Zoom version of (c).}
  \label{fig.maevd.square2o5}
\end{figure}
\subsection{On the solution of problem (\ref{eq.maevd})}
From Sections \ref{sec.intro} and \ref{sec.formulation}, the problem (\ref{eq.maevd}) we consider is to find the ground state solution of
\begin{align}
  \begin{cases}
    u\leq0 \mbox{ and convex}, \lambda>0,\\
    \det\D^2u=-\lambda u |u|^{d-1}\mbox{ in } \Omega,\\
    u=0 \mbox{ on } \partial\Omega,\\
    \displaystyle\int_{\Omega}|u|^{d+1}d\x=1.
  \end{cases}
  \label{eq.maevd.num}
\end{align}
The existence and uniqueness of solutions to problem (\ref{eq.maevd.num}) have been proved in \cite{lions1985two} and \cite{le2017eigenvalue}, assuming that $\Omega$ is strictly convex, bounded and that $\partial\Omega$ is sufficiently smooth (see \cite{lions1985two}, \cite{le2017eigenvalue} for details). The first problem (\ref{eq.maevd.num}) we consider is the one where $d=2$ and $\Omega=\left\{\x=\{x_1,x_2\},x_1^2+x_2^2<1\right\}$. In this particular case, the solution is \emph{radial} and can be computed very accurately, verifying $u(\mathbf{0})=-1.1585...$ and $\lambda=7.4897...$. When applying the methodology discussed in Sections \ref{sec.timediscretize} and \ref{sec.spacediscretize} to the solution of the above test problem (using triangulations of the disk like the one in Figure \ref{fig.mesh} (right), and $\varepsilon=\tau=h^2$), one obtains the results reported in Table \ref{tab.maevd.disk} and Figure \ref{fig.maevd.disk}.
The results reported in Table \ref{tab.maevd.disk} and Figure \ref{fig.maevd.disk} suggest that $\lambda_h-\lambda=O(h)$ and uniform convergence of $u_h$ to $u$ as $h\rightarrow0$.

\begin{rmrk}
  When applying the discrete analogue of the SQP algorithm (\ref{eq.maevd.SQP.0}), (\ref{eq.maevd.SQP.1}), we took $\|u_{k+1}-u_k\|_{0h}\leq 10^{-10}$ as stopping criterion. The related average number of SQP iterations is 10, typically.
\end{rmrk}

We consider now the numerical solution of problem (\ref{eq.maevd.num}) for $\Omega$ defined by (\ref{eq.square2o5}). The related numerical results (some of them reported in Table \ref{tab.maevd.square2o5} and Figure \ref{fig.maevd.square2o5}) confirms those obtained for the unit disk: we have, in particular, uniform convergence of $u_h$ to $u$, and $\lambda_h\approx \lambda-ch$, with $\lambda\approx 6.4$ and $c\approx26$.

\begin{table}[t]
  \begin{tabular}{|c|c|c|c|c|c|c|c|c|}
    \hline
    $h$ & Iteration \# & $\|u_h^{n+1}-u_h^n\|_{0h}$ & $L^2$-error & rate & $L^{\infty}$-error & rate & $\lambda_h$ & $\min u_h$  \\
    \hline
    1/10 & 355 & $9.96\times 10^{-8}$ & $1.30\times 10^{-1}$ &  & $1.15\times10^{-1}$ & &2.55 & -2.5539 \\
    \hline
    1/20 & 385 & $9.99\times 10^{-8}$ & $9.94\times 10^{-2}$ &  0.39& $1.01\times10^{-1}$ &0.19 &2.92 & -2.6043 \\
    \hline
    1/40 & 819 & $9.99\times 10^{-8}$ & $6.64\times 10^{-2}$ & 0.58 & $6.82\times10^{-2}$ & 0.57& 3.36 & -2.5974 \\
    \hline
    1/80 & 2584 & $9.99\times 10^{-8}$ & $3.71\times 10^{-2}$ & 0.84 & $3.91\times10^{-2}$ &0.80 & 3.56 & -2.6081 \\
    \hline
  \end{tabular}
  \caption{Problem (\ref{eq.mabg.num}) with $\Omega=\left\{\x=\{x_1,x_2\},x_1^2+x_2^2<1\right\}$ and $C=10.5$. Variations with $h$ of the number of iterations necessary to achieve convergence (2nd column), of the $L^2$ and $L^{\infty}$ approximation errors and of the associated convergence rates (column 4, 5, 6 and 7), of the computed eigenvalue (8th column) and of the minimal value of $u_h$ over $\Omega$ (that is $u_h(\mathbf{0})$) (9th column). At $C=10.5$, the radial solution verifies $\lambda=3.76...$ and $u(\mathbf{0})=-2.628...$.}
  \label{tab.mabg.disk}
\end{table}

\begin{figure}[t]
  (a)\includegraphics[width=0.35\textwidth]{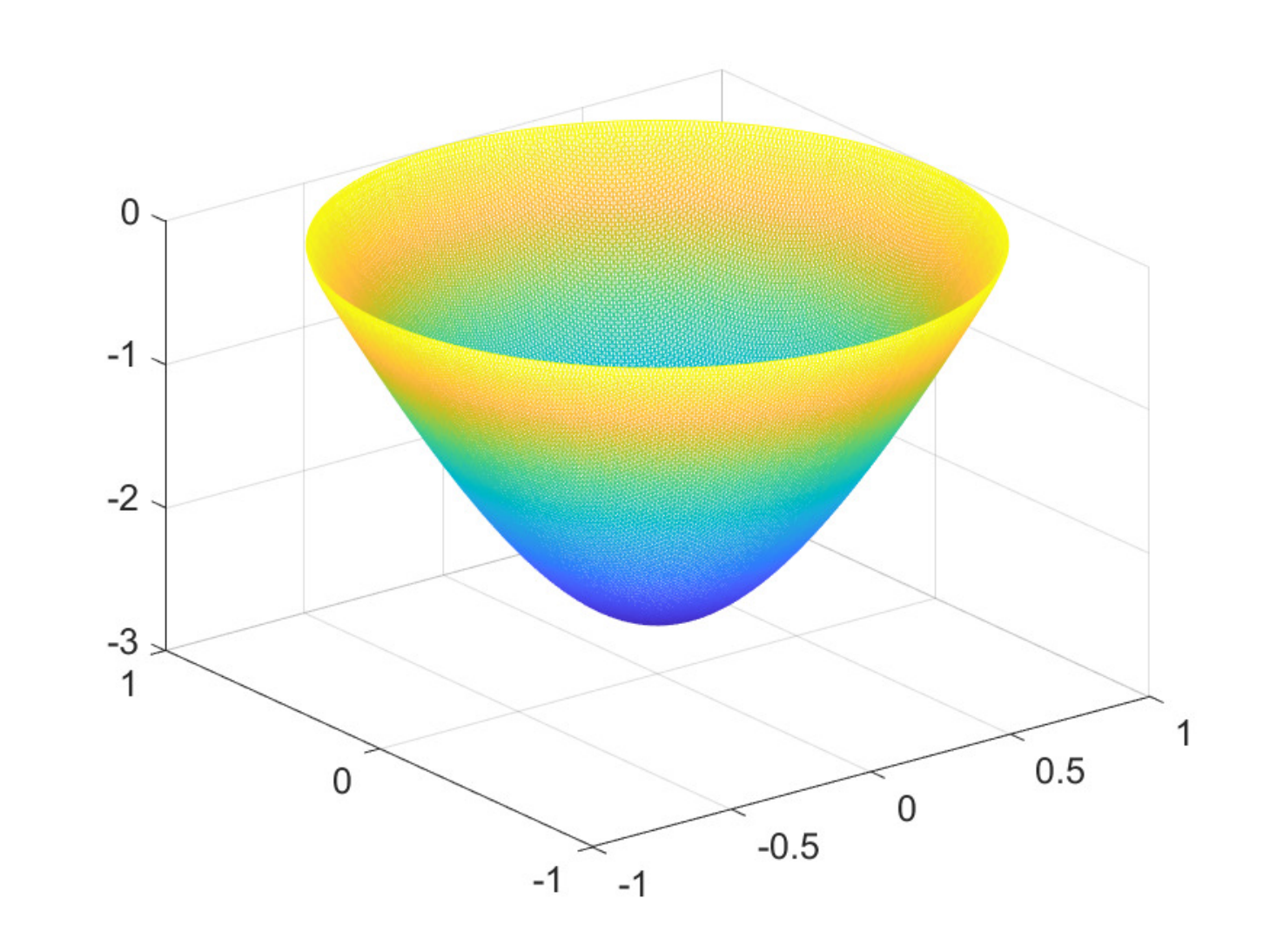}
  (b)\includegraphics[width=0.35\textwidth]{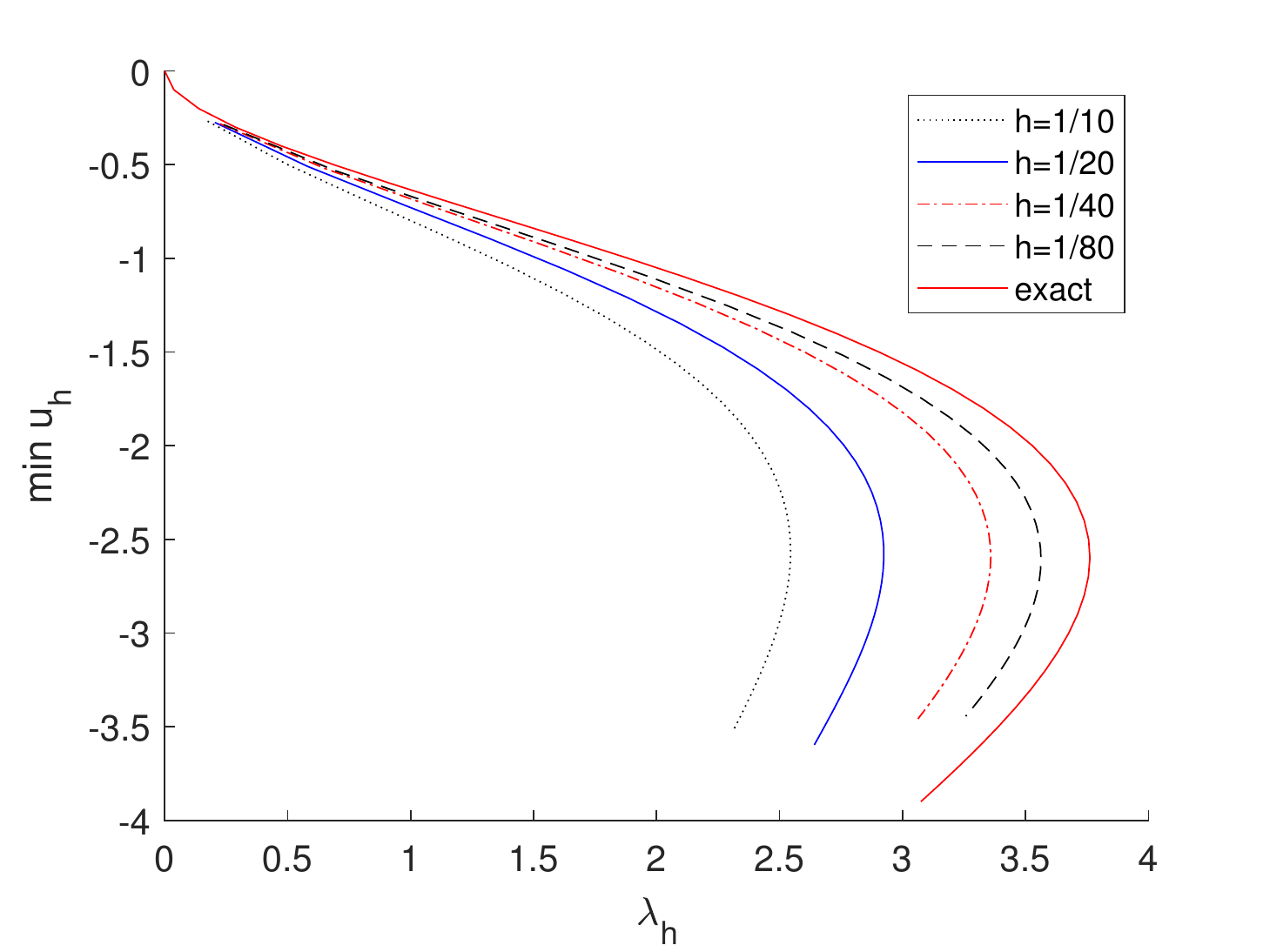}\\
  (c)\includegraphics[width=0.35\textwidth]{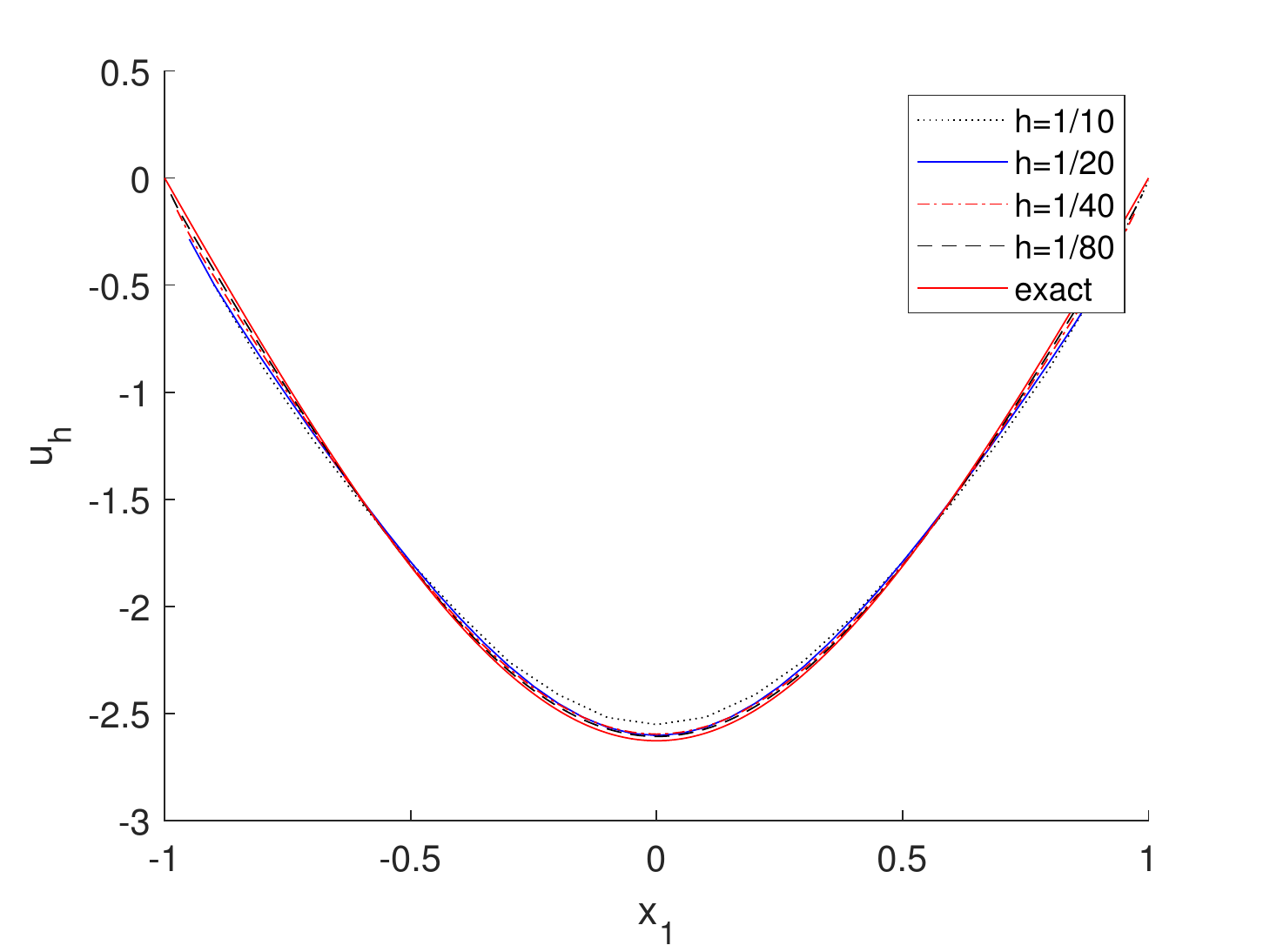}
  (d)\includegraphics[width=0.35\textwidth]{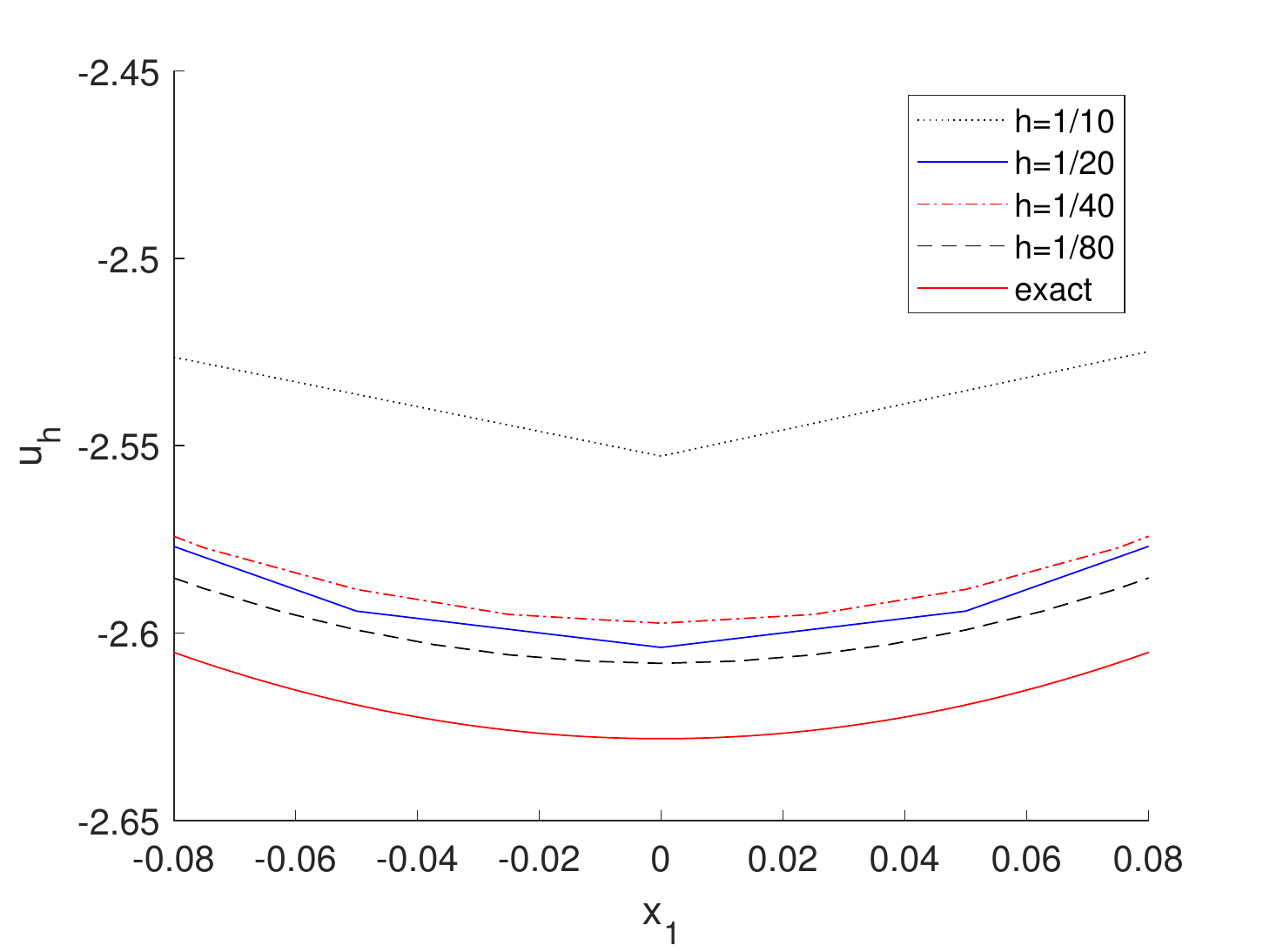}
  \caption{Problem (\ref{eq.mabg.num}) with $\Omega=\left\{ \x=\{x_1,x_2\},x_1^2+x_2^2<1\right\}$. (a) Graph of the approximate solution $u_h$ computed with $h=1/80$ with $C=10.5$. (b) Bifurcation diagrams of the exact solution and of the approximated solutions for $h=1/10, 1/20, 1/40$ and $1/80$. (c) Graphs of the restrictions of $u_h$ and $u$ to the line $x_2=0$ for $h=1/10,1/20,1/40$ and $1/80$. (d) Zoom version of (c)}
  \label{fig.mabg.disk}
\end{figure}

\begin{table}
  \begin{tabular}{|c|c|c|c|c|}
    \hline
    $h$ & Iteration \# & $\|u_h^{n+1}-u_h^n\|_{0h}$ & $\lambda_h$ & $\min u_h$  \\
    \hline
    1/20 & 1865 & $9.98\times10^{-8}$ & 33.35 & -2.5503 \\
    \hline
    1/40 & 1759 & $9.98\times10^{-8}$ & 37.83 & -2.5842 \\
    \hline
    1/80 & 5466 & $9.99\times10^{-8}$ & 41.28 & -2.5990 \\
    \hline
  \end{tabular}
 \caption{Problem (\ref{eq.mabg.num}) with $\Omega$ defined by (\ref{eq.square3}). Variations with $h$ of the number of iterations necessary to achieve convergence (2nd column), of the computed eigenvalue (4th column) and of the minimal value of $u_h$ over $\Omega$ (that is $u_h(\mathbf{0})$) (5th column).}
  \label{tab.mabg.square3}
\end{table}

\begin{figure}
  (a)\includegraphics[width=0.35\textwidth]{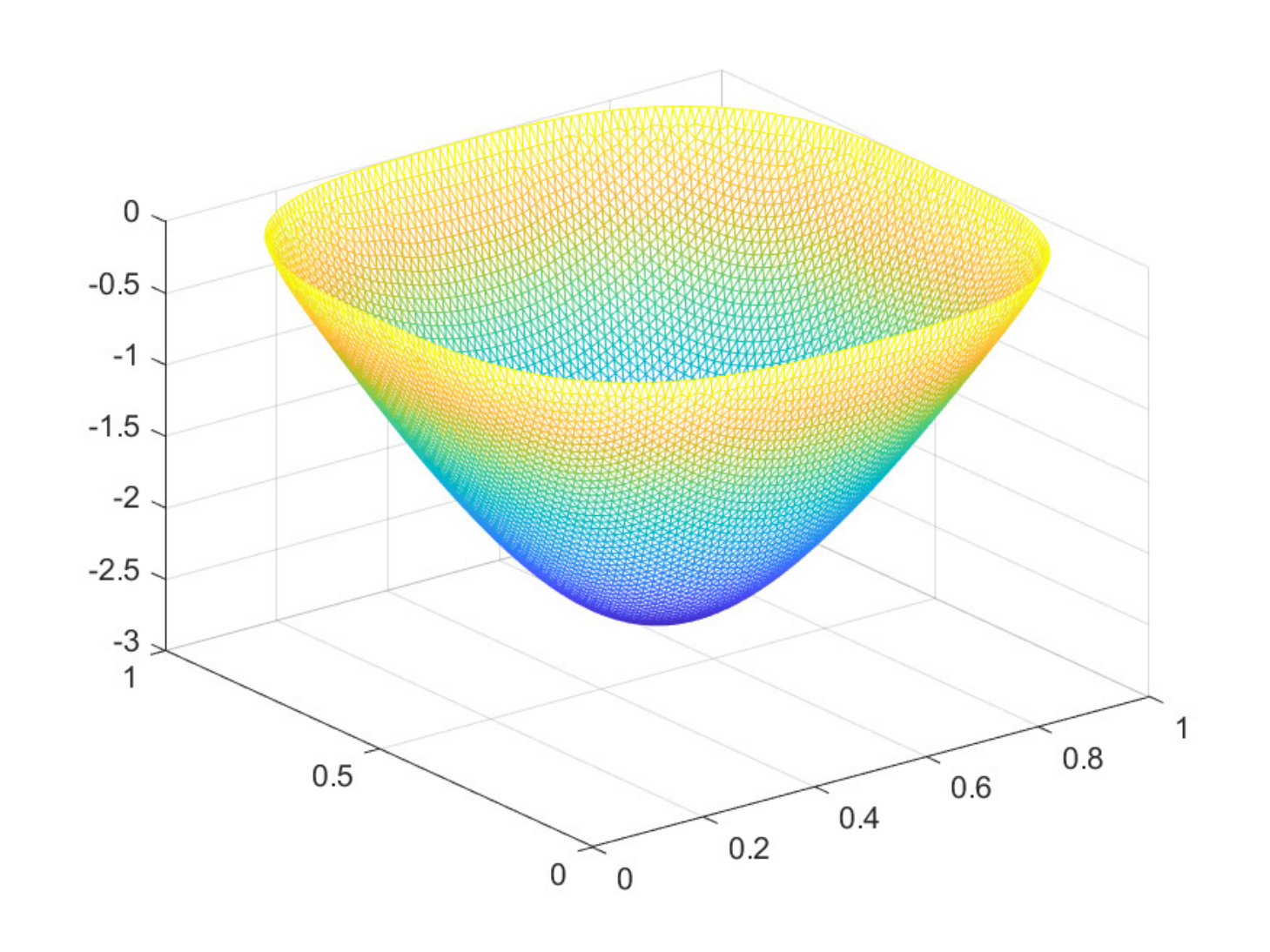}
  (b)\includegraphics[width=0.35\textwidth]{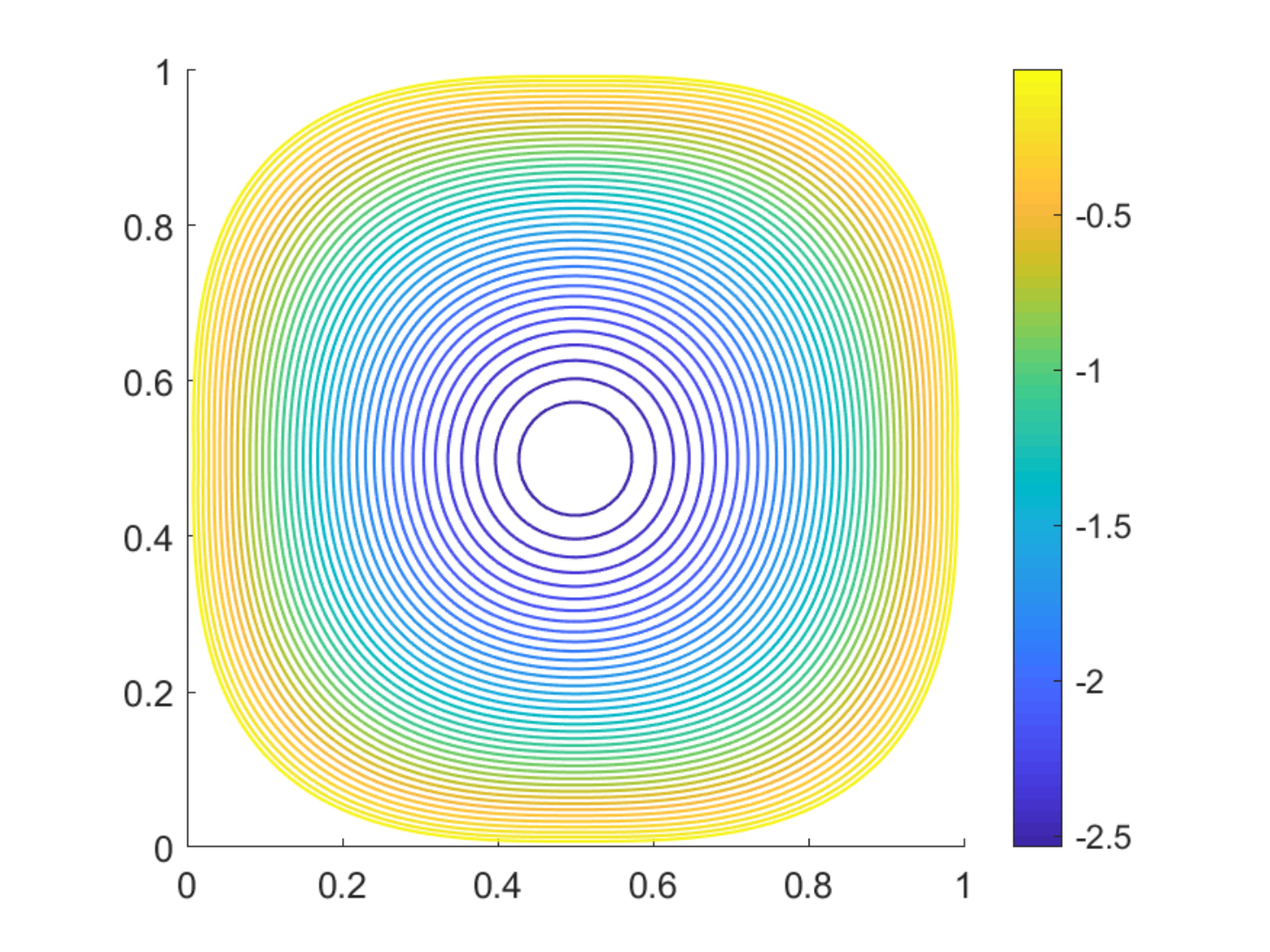}
  (c)\includegraphics[width=0.35\textwidth]{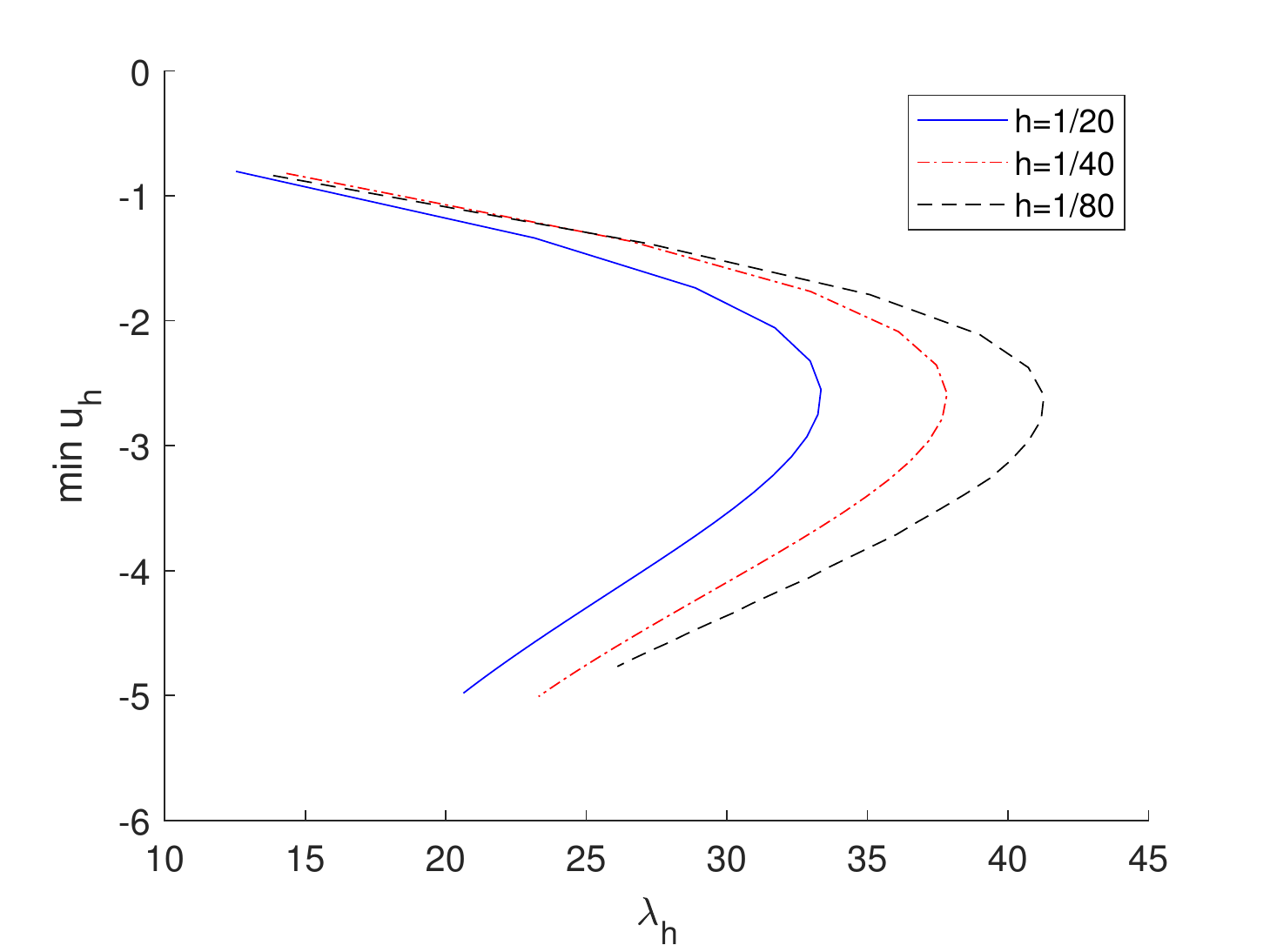}
  (d)\includegraphics[width=0.35\textwidth]{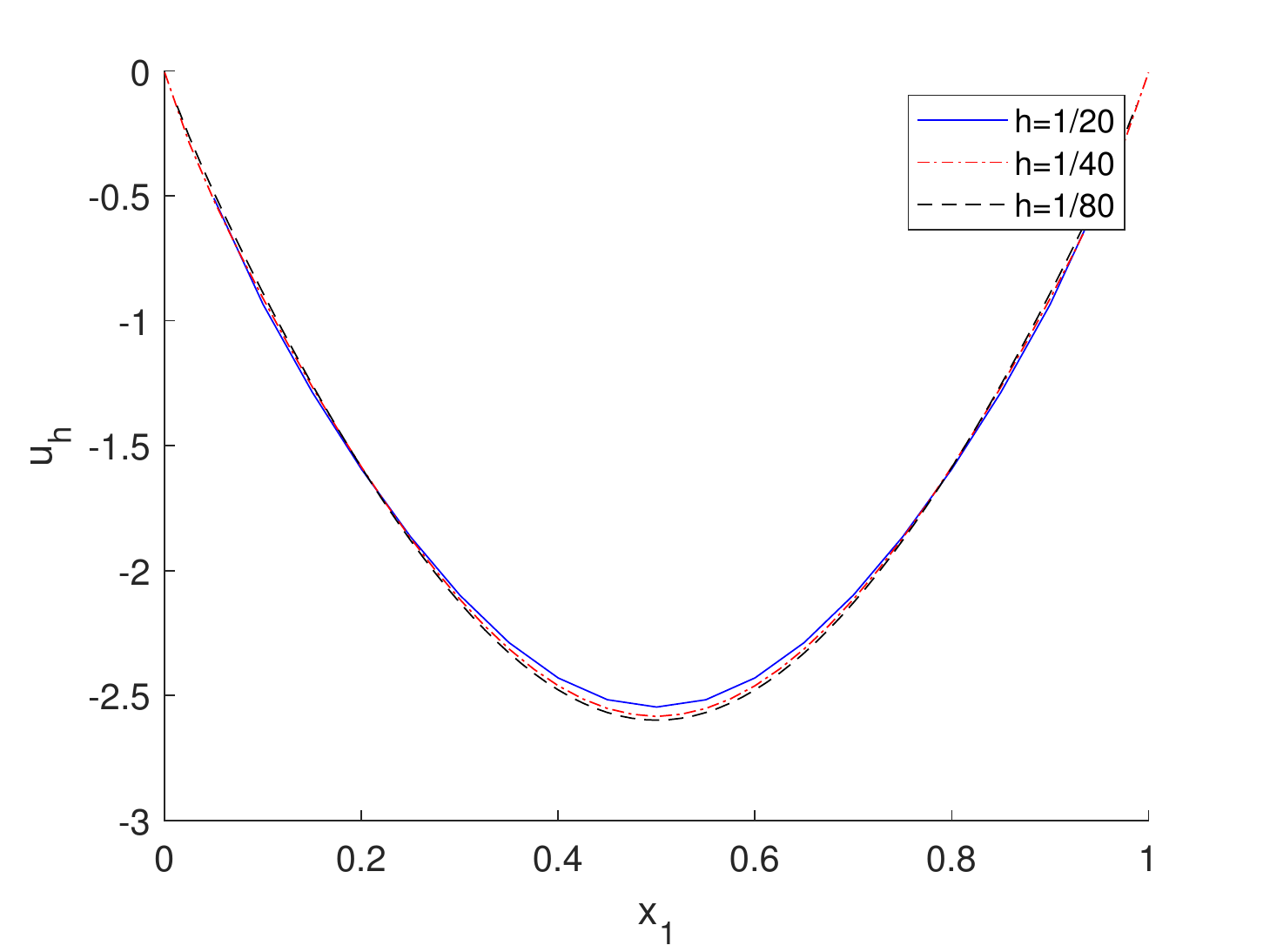}
  \caption{Problem (\ref{eq.mabg.num}) with $\Omega$ defined by (\ref{eq.square3}). (a) Graph of the approximate solution $u_h$ computed with $h=1/80$ for $C=3$. (b) Contours of the approximate solution $u_h$ computed with $h=1/80$ for $C=3$. (c) Bifurcation diagrams of the approximated solutions for $h=1/20, 1/40$ and $1/80$. (d) Graphs of the restrictions of $u_h$ to the line $x_2=1/2$ for $h=1/20,1/40$ and $1/80$.}
  \label{fig.mabg.square3}
\end{figure}
\subsection{On the solution of problem (\ref{eq.mabg}), (\ref{eq.mabg.sys})}
To conclude our investigations, we are going to present some of the results we obtained when applying the methodology we discussed in Sections \ref{sec.timediscretize} and \ref{sec.spacediscretize} to the numerical solution of problem (\ref{eq.mabg}), (\ref{eq.mabg.sys}), namely
\begin{align}
  \begin{cases}
    u\leq0, \lambda>0,\\
    \det\D^2u=\lambda e^{-u} \mbox{ in } \Omega,\\
    u=0 \mbox{ on } \partial\Omega,\\
    \displaystyle\int_{\Omega}\left( e^{-u}-1\right)d\x=C (>0).
  \end{cases}
  \label{eq.mabg.num}
\end{align}
Before presenting the results we obtained when solving various problems of (\ref{eq.mabg.num}) type, let us observe that for a given $\Omega$, the related problems (\ref{eq.mabg.num}) form a family, parameterized by $C$, of nonlinear eigenvalue problems. As done in related situations, we are going to take an \emph{incremental} approach, where one considers a sequence $(C_q)_{q\geq0}$ defined by
$$
\begin{cases}
  C_0=0,\\
  C_{q+1}=C_q+\Delta C,
\end{cases}
$$
with $\Delta C>0$ and 'small'. Next, we solve problem (\ref{eq.mabg.num}) for $C=C_q$ using the method discussed in Sections \ref{sec.timediscretize} and \ref{sec.spacediscretize}, and denote its solution by $(u_{C_q},\p_{C_q},\lambda_{C_q})$. When solving (\ref{eq.mabg.num}) for $C=C_{q+1}$ with $q\geq1$, we advocate initializing our time-stepping method with $(u^0,\p^0)=(u_{C_q},\p_{C_q})$. If $q=0$, we advocate taking $(u^0,\p^0)=\left(-\frac{\Delta C}{|\Omega|},\mathbf{0}\right)$.

The first test problem we are going to discuss is (as expected) the particular case of (\ref{eq.mabg.num}) where $\Omega$ is the unit disk of $\RR^2$ centered at $(0,0)$, that is
\begin{align}
  \Omega=\left\{\x=\{x_1,x_2\},x_1^2+x_2^2<1\right\}.
  \label{eq.disk}
\end{align}
The \emph{bifurcation diagram} of the \emph{radial solution} of problem (\ref{eq.mabg.num}), (\ref{eq.disk}) has been visualized in Figure \ref{fig.mabg.bifurcation}, the \emph{turning point} corresponding to $u(\mathbf{0})=-2.5950...,\lambda=3.7617...$ and $C=10.228...$. The computational methodology discussed in Sections \ref{sec.timediscretize} and \ref{sec.spacediscretize} has been applied to the solution of problem (\ref{eq.mabg.num}), (\ref{eq.disk}), using: (i) $\varepsilon=h^2$. (ii) $\tau=h^2/4$ if $h=1/10$, $\tau=h^2$ if $h=1/20, 1/40, 1/80$. (iii) $\Delta C=0.5$. (iv) Unstructured isotropic finite element triangulations like the one in Figure \ref{fig.mesh} (right). (v) $\left\|u_h^{n+1}-u_h^n\right\|_{0h}<10^{-7}$ as stopping criterion of our operator-splitting based time-stepping method. Related numerical results are  reported in Table \ref{tab.mabg.disk} and Figure \ref{fig.mabg.disk}.
The results reported in Table \ref{tab.mabg.disk} and Figure \ref{fig.mabg.disk} suggest that $\lambda-\lambda_h\approx 16h$ and that $\lim_{h\rightarrow0} u_h=u$, uniformly. They suggest also that the values of $u_h(\mathbf{0})$ corresponding to the turning point do not vary much with $h$ (see Figure \ref{fig.mabg.disk}(b)).

The second (and last) test problem we consider is problem (\ref{eq.mabg.num}) with
\begin{align}
  \Omega=\left\{ \x=\{x_1,x_2\},|x_1-1/2|^3+|x_2-1/2|^3<1/8\right\}.
  \label{eq.square3}
\end{align}
When applying the computational methodology discussed in Sections \ref{sec.timediscretize} and \ref{sec.spacediscretize} to the solution of problem (\ref{eq.mabg.num}), (\ref{eq.square3}), we used: (i) $\varepsilon=h^2$. (ii) $\tau=h^2/32$ if $h=1/20$, $\tau=h^2/8$ if $h=1/40$ and $1/80$. (iii) $\Delta C=0.5$. (iv) Unstructured isotropic finite element triangulations qualitatively like the one in Figure \ref{fig.square2o5}. (v) $\left\| u_h^{n+1}-u_h^n \right\|_{0h}<10^{-7}$ as stopping criterion of our operator-splitting based time-stepping method. Related numerical results are reported in Table \ref{tab.mabg.square3} and Figure \ref{fig.mabg.square3}. The numerical results suggest uniform convergence of $u_h$ to a limit $u$.

\section*{Acknowledgment}
R. Glowinski acknowledges the support of the Hong Kong based Kennedy Wong Foundation. J. Qian is partially supported by NSF grants. S. Leung is supported by the Hong Kong RGC under Grant 16302819.
%%-----------------------------,
\bibliographystyle{abbrv}
\bibliography{reference}
%%-----------------------------
\end{document}